\documentclass[11pt]{article}\input{amssym.def}\input{amssym}
\usepackage[active]{srcltx}
\usepackage{color}

\usepackage{hyperref}
\hypersetup{
    colorlinks=true, 
    linktoc=all,     
    linkcolor=blue,  
}

\newenvironment{annotacia}{\centerline{\sc Abstract}\vspace{2mm}\narrower\narrower\sf}

\makeatletter\@addtoreset{equation}{section}\makeatother

\def\nn{\nonumber}\def\lb{\label}\def\Mt{M \raisebox{1mm}{$\intercal$}}
\def\be{\begin{equation}}\def\ee{\end{equation}}\def\ba{\begin{eqnarray}}\def\ea{\end{eqnarray}}
\def\tr{{\rm Tr}\,}
\def\Tr#1{{\rm Tr}_{\! R^{\mbox{\scriptsize$(#1)$}}}}
\def\TR#1#2{{\rm Tr}_{\! #2^{\mbox{\,\scriptsize$(#1)$}}}}\def\str#1{\rule[#1mm]{0pt}{1mm}}

\newcounter{theorem}\makeatletter
\@addtoreset{theorem}{section}\makeatother

\newtheorem{prop}[theorem]{Proposition}\newtheorem{rem}[theorem]{Remark}\newtheorem{lem}[theorem]{Lemma}
\newtheorem{def-lem}[theorem]{Definition-Lemma}\newtheorem{def-prop}[theorem]{Definition-Proposition}
\newtheorem{defin}[theorem]{Definition}\newtheorem{theor}[theorem]{Theorem}
\newtheorem{cor}[theorem]{Corollary}

\begin{document}

\title{ }
\begin{center}
{\Large \textbf{	Cayley--Hamilton Theorem for Orthogonal \\[5
pt] Quantum Matrix Algebras}}
	
\vspace{1cm} {\large \textbf{Oleg Ogievetsky$^{\diamond\,\dag\,
}$ and Pavel Pyatov$^{\ast\,\star
}$}}
	
\vskip .8cm $^{\diamond}$Aix Marseille Universit\'{e}, Universit\'{e} de
Toulon, CNRS, \\ CPT UMR 7332, 13288, Marseille, France
	
\vskip .3cm $^{\dag}${I.E.Tamm Department of Theoretical Physics, P.N. Lebedev Physical Institute, Leninsky prospekt 53, 119991 Moscow, Russia}
	
\vskip .3cm $^{\ast}${National Research University "Higher School of Economics", 20 Myasnitskaya street, Moscow
101000, Russia}
	
\vskip .3cm $^{\star}${Bogoliubov Laboratory of Theoretical Physics, JINR, 141980 Dubna, Moscow region, Russia}
\end{center}

\begin{annotacia} 
\noindent For a family of the orthogonal $O(k)$ type Quantum Matrix algebras we establish an
analogue of the Cayley--Hamilton theorem. The form of the Cayley-Hamilton identity is different in 
three cases. First, the cases of odd ($k=2\ell -1$) and even ($k=2\ell$)  heights are different. 
Second, for even height orthogonal 
QM-algebra we derive two versions of the Cayley--Hamilton theorem, one for its positive component  
$O^+(2\ell)$ and another one for the negative component $O^-(2\ell)$. 
In each case we introduce the spectral parameterization of the coefficients of the  
Cayley--Hamilton identity by the `eigenvalues' of the quantum matrices.
\end{annotacia}

\newpage
\tableofcontents
\bigskip\bigskip\bigskip

\section{Introduction}\lb{sec1}
In the present work we continue investigation of the quantum matrix (QM-) algebras --- the unital 
associative $\Bbb C$-algebras generated by the components of a square matrix $M$ subject to 
specific quadratic defining relations which are determined by the Yang--Baxter matrix (also called 
$R$-matrix ) 
representations of the braid groups. These algebras were introduced in late 90-ies in \cite{Hl,IOP1}
with an aim to develop a unified approach to the algebras of quantized functions on the Lie groups 
(the so-called FRT-algebras, see \cite{FRT}) and the algebras of quantized vector 
fields over the Lie groups
(the so-called Reflection Equation (RE-) algebras, see \cite{KS}).

\vskip .1cm
Main structural properties of the QM-algebra depend on the type of the 
Yang--Baxter matrix  used in its definition. The Yang-Baxter matrices realizing 
representations of the Iwahori--Hecke  algebras \cite{Iw,K} define a family of the Hecke type QM-algebras. Typical representatives of this family are
the QM-algebras describing the standard general linear quantum (super-)groups.  
Correspondingly, one can further subdivide the Hecke type QM-algebras into  
subfamilies of $GL(m|n)$ type, $m, n \geq 0$.\footnote{We do not know any Hecke type QM-algebra 
that do not fall into this classification.} 

The QM-algebras of the $GL(m)$ type ($n=0$ case) are 
the most known and well investigated. General structure results for them, including the 
Cayley--Hamilton identities, the Newton and Wronski relations in their  characteristic subalgebras, 
and the Cayley--Hamilton--Newton matrix identities unifying  both previously mentioned results were 
obtained in \cite{PS,GPS,IOP,IOP1,IOP2}. Many particular examples of these type QM-algebras 
were studied in \cite{Gr,BCC,Gou,GZB,EOW,NT,Zh,JW,OV,OV1}.

For the QM-algebras of the type $GL(m\vert n)$, $m,n>0$,
the general structure theory was 
developed in \cite{GPS2,GPS3}. The results for the classical supermatrices were obtained earlier in 
\cite{KT1,KT2} and for the universal enveloping algebras of the Lie superalgebras in \cite{JGr,IWG}.
 
\vskip .1cm
In the present paper we study the QM-algebras associated with the Birman--Murakami--Wenzl 
(BMW, for short) algebras \cite{BW,M1} – finite dimensional quotient algebras of the 
group algebras of the braid groups  which may be treated as deformations of the 
Brauer algebras \cite{Br}. This 
research was initiated in the unpublished paper \cite{OP}. Later, in a series of papers 
\cite{OP-BMW,OP-SpCH,OP-reciprocal} we 
developed the general structure theory of the BMW type QM-algebras, including the notion 
of its characteristic 
subalgebra and a generalization of the quantum matrix $M$ powers. 
These notions are the main ingredients 
in a construction of the characteristic or, Cayley-Hamilton identities for the QM-algebras. Namely, 
the characteristic 
identity is a linear combination of the generalized matrix powers of the quantum matrix $M$, 
while the scalar coefficients of this combination are the elements of the characteristic subalgebra.

n \cite{OP-SpCH} we separated two subfamilies of the BMW type QM-algebras ---the algebras of the 
orthogonal $O(k)$ and the symplectic $Sp(2k)$ types. For the symplectic QM-algebras we  
succeeded in deriving the Cayley-Hamilton identities. The case of 
the orthogonal QM-algebras turns out 
to be more sophisticated: the generators of the characteristic subalgebra in this case satisfy 
additional, so-called reciprocal relations which we have obtained in \cite{OP-reciprocal}. 

In the present work 
we derive the Cayley–Hamilton identities for
the orthogonal QM-algebras. We also construct 
spectral parameterizations of their characteristic subalgebras. 
By the spectral parameterization we mean the homomorphism from the characteristic 
subalgebra to the commutative polynomial algebra of the so-called `spectral variables'. These 
homomorphisms are governed by the nature of the Cayley--Hamilton identities; the
Cayley--Hamilton identities factorize after the introduction of the spectral variables, 
 which allows them to be interpreted
as the eigenvalues of the quantum matrix $M$. Thus we complete the 
general structure theory for the orthogonal and symplectic types 
QM-algebras. The particular 
examples were considered earlier in \cite{BGr,Gr,Gou,Mudr}. 

\vskip .2cm 
Since this work, in a sense, summarizes an entire piece of our research 
let us briefly outline some further research
directions before moving on to the detailed presentation of the paper.

First of all, it is worthwhile
to discuss low-dimensional 
examples of the characteristic identities in the FRT- and RE- algebras in detail. We 
plan to address this in \cite{OP-examples}. 
In contrast to the $GL$ case, the Poincar\'e
series for the QM-algebras in the BMW case do not coincide with the classical ones 
and it would be interesting to trace the general pattern. 

Next, it is known  \cite{GPS3} that in the Hecke case the characteristic subalgebra 
admits a basis in which the multiplication 
reproduces the Littlewood--Richardson rule for the Schur symmetric functions. It would 
be interesting to understand an analogue of this phenomenon in the BMW case.

\vskip .1cm
Finally, in addition to families of the orthogonal and symplectic QM-algebras, the BMW type contains 
orthosymplectic QM-algebras, unifying orthogonal and symplectic
types. Their corresponding Yang-Baxter matrices were described in \cite{BSh,I}. This case appears labor-intensive to us, since it must combine characteristic 
identities for $q$-symplectic and $q$-orthogonal Lie algebras, which look completely different. We hope to 
return to this 
$q$-orthosymplectic concordance soon. This will exhaust the quantization of 
`classical' Lie (super)-algebras having a `vector' representation. 
Some particular orthosympectic examples were considered in \cite{JGr,GrJ,Gou2,GI}. 

\medskip
The paper is organized as follows.

\vskip .2cm
Section \ref{sec2} contains necessary information about Yang--Baxter matrices and Yang--Baxter 
representations of the braid groups and, more specifically, of the Birman--Murakami--Wenzl
(BMW-, for short) algebras -- certain finite-dimensional quotient algebras of the group algebras of the 
braid groups. Namely, we introduce the notion of a strict skew-invertible
Yang--Baxter matrix $R$ and the corresponding $R$-trace. We briefly recall the definitions 
and properties of the finite height BMW-type Yang--Baxter matrices  
and the subfamily of the $O(k)$-type (that is, height $k$ orthogonal type) 
Yang--Baxter matrices. Further we recall the definitions of the compatible Yang--Baxter pair 
$\{ R,F\}$ and the twisted Yang--Baxter matrix. All these concepts will be used for the definitions 
of the quantum matrix algebras and related objects. 

\vskip .2cm
The quantum matrix (QM-, for short) algebra ${\cal M}(R,F)$ associated to a compatible
Yang--Baxter pair $\{ R,F\}$ is introduced in section \ref{sec4}. 

\vskip .1cm
In subsection \ref{subsec4.1}
we specialize to the case of the orthogonal QM-algebra of the type $O(k)$, and define its 
certain commutative subalgebra ${\cal C}(R,F)$. The subalgebra ${\cal C}(R,F)$ is called
`characteristic' because the coefficients 
of the future characteristic (Cayley--Hamilton) identity for the matrix $M$ of generators  
of the QM-algebra take values in it. We introduce three sets of generators of the 
subalgebra ${\cal C}(R,F)$ containing the sequences of the power sums, of the complete 
sums and of the 
elementary sums. In addition, each set contains the element $g$ called contraction. It is exactly the 
appearance of the element $g$ which makes a difference between the BMW and the Hecke cases. 
In the orthogonal case the sequence of elementary sums is finite. Moreover the 
elementary sums, together with $g$, become dependent and satisfy the quadratic (in elementary 
sums) reciprocal
relations, derived in \cite{OP-reciprocal}. 

\vskip .1cm
In subsection \ref{subsec4.2} we work with the set 
of matrices $M^{\alpha^j}$ 
indexed by elements
$\alpha_j$ of the braid groups ${\cal B}_j$ for all $j=0,1,2...$. Their components  belong to the
QM-algebra ${\cal M}(R,F)$. We introduce a special `$\star$-multiplication' for matrices of our set.
We extend the linear span of our set 
by supplying it with a multiplication by the elements of the 
algebra ${\cal C}(R,F)$. The result is an algebra ${\cal P}(R,F)$ which we call the quantum matrix 
powers algebra.
This algebra is associative and, in the BMW case, commutative. We discuss the set of generators of 
this algebra in the BMW case and define the `$\star$-ìnverse' of the quantum matrix $M$.

\vskip .1cm
In subsection \ref{components} we remind the reader of the resolution of the 
reciprocal relations from which becomes possible for an invertible $g$.
 
\vskip .2cm
Section \ref{sec5} contains the main results of the paper. 
We construct the generalization of the Cayley--Hamilton identity for the orthogonal type quantum matrices.

\vskip .1cm
Subsection \ref{subsec5.1} contains some preparatory results. We introduce two 
sequences of elements from 
${\cal P}(R,F)$. Roughly speaking these sequences are the compositions of 
the usual and wedge matrix powers of the quantum matrix $M$. We establish 
recursive relations for them 
and present the solutions of these recursions. Proposition \ref{prop5.2} prepares then the main
ingredient for the the characteristic identities, containing exclusively the matrix powers of
of the quantum matrix $M$. 

\vskip .1cm
The construction of the Cayley--Hamilton identity splits into three different cases. 

\vskip .1cm
In subsection \ref{subsec5.2} we  investigate the even height case, $k=2\ell$.
We separate two components of the orthogonal QM-algebra of the type $O(2\ell)$. This 
separation is related to the decomoposition of the classical group $O_{2\ell}$ into the union of 
the two components  $O_{2\ell}^+$ and  $O_{2\ell}^-$. The Cayley--Hamilton 
identities look different for the components, they are given in theorem \ref{theorem5.6}. 

\vskip .1cm
In subsection \ref{subsec5.3} we present the Cayley--Hamilton identity for the orthogonal
QM-algebra $O(2\ell-1)$. It is given in theorem \ref{theorem5.8}.

\vskip .2cm
Section \ref{sec6} introduces a set of variables, referred to as ‘spectral variables’, intended to 
represent the eigenvalues of the 
quantum matrix. These variables are defined in a consistent with the reciprocal relations manner,
they satisfy the quadratic relations. 

\vskip .1cm
In subsection \ref{subsec6.1} we introduce three  
homomorphic maps from the characteristic subalgebras 
of the QM-algebras ${O}_{k}(R,F)$ into the algebras of polynomials in the spectral variables. 
These maps are defined differently in the three cases considered in subsections \ref{subsec5.2} and 
 \ref{subsec5.3}.
The maps allow to express, in each case, the elementary sums ${\goth e}_i$ in terms of  the spectral 
variables and
therefore we refer to these homomorphic maps as spectral parameterizations.
Using these maps we define central extensions  of the QM-powers algebras ${\cal P}_{k}(R,F)$. 
Then the Cayley--Hamilton identities from theorems \ref{theorem5.6} and 
\ref{theorem5.8}, considered as 
relations in the central extensions of the QM-powers algebras, acquire completely factorized forms. 

\vskip .1cm
In subsection \ref{subsec6.2} we derive the 
spectral parameterizations for the other two generating sets 
of the characteristic subalgebra -- the power sums ${\goth p}_i$ and the complete sums 
${\goth h}_i$. 

\vskip .1cm
The spectral parameterization is the second main result of the present work.

\vskip .1cm
As in the previous paper \cite{OP-reciprocal}, we shall call a matrix solution of the Yang--Baxter
equation not as it is commonly named `R-matrix' but `Yang--Baxter matrix', to avoid ugly   
combinations like `R-matrix R'. Also, a braid group representation, built on a 
Yang--Baxter matrix, we call a `Yang--Baxter' representation.

\vskip .1cm
The symbol $\blacktriangleright$ denotes ends of remarks as well as ends of proofs of auxiliary 
lemmas placed inside the proofs of theorems (whose ends are 
denoted by $\blacksquare$).

\section{Yang--Baxter matrices}\lb{sec2}
 In this section we remind some basic facts about Yang--Baxter representations of the braid group
 of type A. These data will be used in the definition and investigation of the quantum matrix algebras. 
 A detailed presentation of the material can be found in \cite{OP-BMW}, sections 2, 3, and 
 \cite{OP-reciprocal}, section 3.

 \medskip
Let $V$ be a finite dimensional $\Bbb C$-linear space, $\dim V =\mbox{\sc n}$. We fix some basis 
$\{v_i\}_{i=1}^{\rm N}$ in $V$ and the induced bases 
$\{v_{i_1}\otimes\dots\otimes v_{i_n}\}_{i_1,\dots,i_n=1}^{\rm N} $ in spaces 
$V^{\otimes n}$, $n=1,2, \dots$ .

 \medskip
Denote by ${\rm Id}_{_V}$ the identity operator in space $V$. 
For any operator $F\in{\rm End}(V^{\otimes j})$, $j\geq 1$, and for $n\geq j$ denote by $F_i$, 
$1\leq i\leq n-j+1$, the operators 
${\rm Id}_{_V}^{\,\otimes (i-1)}\otimes F\otimes {\rm Id}_{_V}^{\,\otimes(n-i-j+1)}\in V^{\otimes n}\,$.
\footnote{Strictly speaking we should write something like $F_{i,n}$ and not simply $F_i$.
However the value of $n$ will be clear from the context and we omit it.  
This remark applies as well to all other objects acting in the tensor product of multiple copies 
of the space $V$.}
\begin{defin}
An operator $F\in{\rm Aut}(V^{\otimes 2} )$ is called Yang--Baxter matrix if the pair of the operators $F_1$, $F_2$ fulfills the braid relation
\be\lb{braidX}
F_1 F_2 F_1 = F_2 F_1 F_2.
\ee
\end{defin}
An important example of the Yang--Baxter matrix is the {\em flip} $P$:
\be
\lb{P}
P(u\otimes v) := v\otimes u \;\; \forall\, u,v\in V.
\ee
\vskip -.2cm
\begin{defin}
The Yang--Baxter matrix $F$ is called  skew invertible if there exists a unique operator
$D_F\in {\rm End}(V)$ such that 
\be
\lb{Tr2R}
\tr_{\!(2)} (D_F)_2\, F_{1} = ({\rm Id_{_V}})_1
\ee
holds true. 
\end{defin}
Here the symbol ${\rm Tr}_{(i)}$ stands for the trace over an $i$-th component space in the tensor power of the spaces $V$. 
\begin{defin}\label{sski}
The Yang--Baxter matrix $F$ is called strict skew invertible if it is skew invertible and the operator
$D_F$ is invetrible. In that case the inverse Yang--Baxter matrix $F^{-1}$ is also strict skew invertible.
\end{defin}
The definitions \ref{P} and \ref{sski} differ from the standard definitions of the (strict) skew invertibility (see, e.g., \cite{Og}) but are equivalent to them.

\begin{defin}	Let $W$ be  a $\Bbb{C}$-linear space.
To a skew invertible Yang--Baxter matrix $F$ one associates a $\Bbb{C}$-linear map
$
{\rm Tr\str{-1.3}}_{\! F}:\; {\rm End}(V)\otimes W\,\rightarrow \, W 
$
called the F-trace:
\be
\lb{R-Tr}
{\rm Tr\str{-1.3}}_{\! F}(X): ={\rm Tr}_{_V} (D_F X)  \;\; \forall\,
X\in{\rm End}(V)\otimes W.
\ee
\end{defin}

To give an example, the  flip $P$ is strict skew invertible, $D_P={\rm Id}_{_V}$ thus, the associated $P$-trace coincides with the usual trace  ${\rm Tr}_{_V}$.
\bigskip

We remind that, in the Artin's presentation, braid group ${\cal B}_n$ is given by a set 
of generators $\{\sigma_i\}_{i=1}^{n-1}$ satisfying defining relations
\ba
\lb{braid group}
\sigma_i \sigma_{i+1} \sigma_i \, =\,  \sigma_{i+1} \sigma_i \sigma_{i+1},&&
\sigma_i \sigma_j \, =\,  \sigma_j \sigma_i\qquad   \forall\; i,j:\; |i-j|>1 ,
\ea

To any Yang--Baxter matrix $F$ one associates representations of the braid groups  
$\rho_F: {\cal B}_n\rightarrow {\rm Aut}(V^{\otimes n})$, $n=2,3, \dots$, which are defined by 
\be
\rho_F(\sigma_i) = F_i, \quad i=1,\dots n-1.
\ee
We call such representations of the braid groups the {\em Yang--Baxter representations}.

 \medskip
The Yang--Baxter representations which realize the  
{\em Birman--Murakami--Wenzl (BMW) quotients} ${\cal W}_n(q,\mu)$ of the group algebras 
${\Bbb C}[{\cal B}_n]$ \cite{BW,M1} are of particular interest to us in the present work. We do not 
give a precise description of the BMW algebras here, we describe only the necessary facts. For the 
definition and proofs, the interested reader is referred  to \cite{OP-BMW}, section 2. The notation 
we use here is also borrowed from that paper.

 \medskip
The BMW algebras depend on two complex parameters $q$ and $\mu$ subject to conditions
\be
\lb{init-cond}
q\notin \{0,\pm 1\},\quad \mu\notin\{0, q, -q^{-1}\}.
\ee
Imposing additional restrictions on the parameters
one can construct two remarkable sets of idempotents in ${\cal W}_n(q,\mu)$: $\{a^{(i)}\}_{i=2}^n$ 
and $\{s^{(i)}\}_{i=2}^n$, called, respectively, {\em $q$-antisymmetrisers} 
and {\em $q$-symmetrisers}. 
Namely, the $q$-antisymmetrizers are defined for
\be
\lb{mu-a}
q^{2j}\,\neq\, 1\, , \quad \mu \neq -q^{-2j+3} \;  \quad \forall\;j = 2,3,\dots , n\ ,
\ee
and the $q$-symmetrizers for
\be
\lb{mu-s}
q^{2j}\,\neq\, 1\, , \quad\;\; \mu \neq q^{2j-3}\;\;\; \quad \forall\;j = 2,3,\dots , n\ .
\ee

The algebra ${\cal W}_2(q,\mu)$ is abelian.  
Its single generator $\sigma_1$ satisfies the cubic relation
$$
(q-\sigma_1)(q^{-1}+\sigma_1)(\mu-\sigma_1)=0,
$$
from which 3 mutually orthogonal idempotents can be constructed 
\be
\lb{3-idemp}
a^{(2)}={(q-\sigma_1)(\mu-\sigma_1)\over (q+q^{-1})(\mu+q^{-1})}, \quad
s^{(2)}={(q^{-1}+\sigma_1)(\mu-\sigma_1)\over (q^{-1}+q)(\mu-q)}, \quad
c^{(2)}={(q-\sigma_1)(q^{-1}+\sigma_1)\over (q-\mu)(q^{-1}+\mu)}.
\ee
This gives rise to the resolution of unity in ${\cal W}_2(q,\mu)$: $1=a^{(2)}+s^{(2)}+c^{(2)}$. 
The third idempotent -- $c^{(2)}$ -- is called {\em the
contractor}. 

 \medskip
The Yang--Baxter matrices which give rise to Yang--Baxter representations of the BMW algebras are called {\em BMW type}. In this text for the BMW type Yang--Baxter matrices we always use notation $R$. Images of the contractor and  of the $q$-antisymmetrizers play major role in considerations below. We use the symbol $K$ for the Yang--Baxter matrix image of the contractor:
\be
\lb{K}
K:= { (q-\mu)(q^{-1}+\mu)\over\mu\,(q-q^{-1})}\, \rho_R(c^{(2)}).
\ee
The chosen normalization of $K$ will be suitable in the sequel. 

 \medskip
If the BMW type Yang--Baxter matrix $R$ is skew invertible, then it is strict skew invertible and the rank of the associated operator $K$ equals 1 \cite{IOP3}: 
\be
\lb{rkK=1}
\mbox{rk}\,K=1.
\ee

We  say that the BMW type Yang--Baxter matrix $R$ is of {\em finite height} $k$  if a 
sequence $\rho_R(a^{(i)})$, $i=1,2,\dots$, terminates at $i=k+1\,$:\,\footnote{This definition can be 
slightly weakened, see \cite{OP-SpCH}, 
section 2.4.  
In general, one does not need the restriction $q^{2(k+1)}\neq 1$ for the parameter $q$.}
\be
\lb{height-k}
\rho_R(a^{(i)})\not\equiv 0\;\; \forall\; i\leq k,\qquad \rho_R(a^{(k+1)})\equiv 0\;\;\Rightarrow \;\;
\rho_R(a^{(k+i)})\equiv 0\;\; \forall\; i> 0.
\ee
The height $k$ condition restricts possible values of the parameter $\mu$ to a pair of essentially 
different possibilities.
The first choice $\mu=-q^{-1-2k}$ gives a family of so called symplectic $Sp(2k)$ type Yang--Baxter 
matrices. The corresponding quantum matrix algebras were investigated in \cite{OP-SpCH}. In the 
present work we are dealing with the second choice.
\begin{defin}\lb{def-O(k)}
A skew invertible BMW type Yang--Baxter matrix $R$ of height $k\geq 2$ is called orthogonal $O(k)$ 
type if the following conditions 
\be
\lb{rank=1}
\mu = q^{1-k}, \qquad
\mbox{\rm rk}\,\rho_R(a^{(k)})=1.
\ee
are satisfied.
\end{defin}

The standard example of the $\,O(k)$ type Yang--Baxter matrix arises from the universal Yang--Baxter 
matrix of the standard $q$-deformation of the orthogonal group \cite{FRT}. This Yang--Baxter matrix 
acts in the tensor square of the $k$-dimensional space and is given by formula
\be
\lb{R-st}
R^{\circ} \! :=\!\!  \sum_{i,j=1}^{k} q^{(\delta_{ij}-\delta_{ij'})}\, E_{ij}\otimes E_{ji} + (q-q^{-1})\!\sum_{i,j=1\atop i<j}^k \left[E_{ii}\otimes E_{jj} -
\,q^{(\rho_j-\rho_i)} E_{ij'}\otimes E_{i'j}\right] ,
\ee	
where $E_{ij}$ are $k\times k$  matrix units,~ $i':=k+1-i$,~~ 
$\rho_i:=-\rho_{i'}:=\frac{k}{2}-i$~  $\forall\;i=1,\dots, \lfloor k/2\rfloor$ and additionally~ $\rho_{(k+1)\over 2}:=0$~ in case if $k$ is odd. 

\medskip
Another concept that we  employ to define a general quantum matrix algebra is the notion of compatible pair of Yang--Baxter matrices. 
\begin{defin}
An ordered pair $\{ R, F\}$ of two Yang--Baxter matrices $R$ and $F$  is
called a compatible Yang--Baxter pair if the following conditions
\be 
R_1\, F_2\, F_1\, =\, F_2\, F_1\, R_2\, ,\qquad R_2\, F_1\, F_2\, =\, F_1\, F_2\, R_1\, 
\lb{sovm}
\ee
are satisfied. The conditions (\ref{sovm}) are called  twist relations. 
\end{defin}

Clearly, $\{ R,P\}$ and $\{ R,R\}$ are compatible Yang--Baxter pairs.

\medskip
Given a compatible Yang--Baxter pair $\{R,F\}$ one can construct a new Yang--Baxter matrix
\be
\lb{R_F }
R_F  := F^{-1} R F\, ,
\ee
called the {\em twisted} Yang--Baxter matrix. The pair $\{ R_F  , F\}$ is again a compatible
Yang--Baxter pair. If $R$ and $F$ are strict skew invertible, then so is $R_F $,
see \cite{OP-BMW}, proposition 3.6.

\section{Quantum matrix algebra in the orthogonal case}\lb{sec4}

In this section we introduce the main objects of our study, the {\em quantum matrix (QM-) algebras} 
of the orthogonal type.
\smallskip 

In subsection \ref{subsec4.1} we start with a general definition of the QM-algebra ${\cal M}(R,F)$ 
associated to a compatible Yang--Baxter pair $\{R,F\}$ of the strict skew invertible Yang--Baxter 
matrices. Depending on the type of the Yang--Baxter matrix $R$ used in the definition, one can 
specialize to QM-algebras of the Hecke and of the BMW types (see \cite{IOP1} and 
\cite{OP-BMW}). In the latter case one can further distinguish subfamilies of the symplectic and the 
orthogonal QM-algebras, and of their super analogues. In \cite{OP-SpCH} we carried out 
investigation of the symplectic $Sp(2k)$ type QM-algebras. Here we continue our study 
in \cite{OP-reciprocal} of the 
orthogonal $O(k)$ type QM-algebras.

 \medskip
A crucial role in the study of the QM-algebra is played by its special commutative subalgebra 
${\cal C}(R,F)$ called the {\em characteristic subalgebra}.  In the end of subsection \ref{subsec4.1} 
we collect information on the characteristic subalgebra of the orthogonal $O(k)$ type QM-algebra 
obtained in \cite{OP-BMW,OP-reciprocal}.

 \medskip
In subsection \ref{subsec4.2} we briefly remind quantum versions of the matrix product and the 
matrix inversion.
This material is based on \cite{OP,OP-BMW}. 

 \medskip
A concluding subsection \ref{components} 
contains results concerning the resolution from \cite{OP-reciprocal} of the reciprocal relations.

\subsection{$O(k)$ type quantum matrix algebra and its characteristic subalgebra}\lb{subsec4.1}

Let $F\in {\rm Aut}(V\otimes V)$ be a skew invertible Yang--Baxter matrix. For any matrix 
$M = \|M_a^b\|_{a,b=1}^{\mbox{\footnotesize\sc n}}$, where $M_a^b$ is an element of some vector space $W$, we introduce a sequence of its {\em matrix 
copies} $M_{\overline i}\in W\otimes{\rm End}(V\otimes V\otimes V\otimes ...)$, $i=1,2,3,\dots\,$:
\be
M_{\overline 1}:=M_1, \quad M_{\overline{i}}:=
F^{\phantom{-1}}_{i-1}M_{\overline{i-1}}F_{i-1}^{-1}\ .
\lb{kopii}
\ee	

\begin{defin}\lb{definitionQMA}
Let $\{R,F\}$ be a compatible pair of the strict skew invertible Yang--Baxter matrices. A general
quantum matrix (QM-)algebra ${\cal M}(R,F)$ is 
given in terms of a set of generators $\{M_a^b\}_{a,b=1,\dots,\mbox{\sc n}}$, subject to defining 
quadratic relations 
\be 
R_i M_{\overline i}M_{\overline{i+1}} = M_{\overline i}M_{\overline{i+1}}R_i, \qquad i=1,2,\dots\, .
\label{qmai}
\ee
The matrix $M$ whose components generate the QM-algebra is called quantum matrix.

The QM-algebra ${\cal M}(R,F)$ is called orthogonal $O(k)$ type if the Yang--Baxter matrix $R$ is orthogonal $O(k)$ type.
\end{defin}

Note that the relations (\ref{qmai}) for different values of index $i$ are pairwise equivalent.
In the BMW case and, in particular, for the orthogonal QM-algebras relations (\ref{qmai}) imply the equalities
\be
\lb{tau2}
K_{i}\, M_{\overline{i}}M_{\overline{i+1}} =
M_{\overline{i}}M_{\overline{i+1}}\, K_{i}\,\, =\,\,\mu^{-2} K_{i}\, g\, 
\qquad 
i= 1,2,\dots\, , 
\ee
where we introduced notation $g$ for quadratic combination of the algebra generators
\be
\lb{tau}
 g:=  \Tr{1,2} \Bigl( M_{\overline{1}}M_{\overline{2}}\, \rho_R(c^{(2)})\Bigr) \, =
 \, { \mu\,(q-q^{-1})\over (q-\mu)(q^{-1}+\mu) }\,
	\Tr{1,2} \left( M_{\overline{1}}M_{\overline{2}}\, K_1\right). 
\ee
It is called  a  {\em contraction} of the quantum matrix $M$. 

\medskip
Here and below the symbol $\Tr{i_1,i_2,...}$ stands for the $R$-trace in the $i_1$-th, $i_2$-th, ...
copies of the space $V$.

 \medskip
Consider a subspace ${\cal C}(R,F)$ in the QM-algebra which is
spanned linearly by the unity and elements
\be
\lb{char}
{\rm ch}(\alpha^{(n)}) := \Tr{1,\dots ,n}(M_{\overline 1}\dots M_{\overline n}\,
\rho_R(\alpha^{(n)}))\ ,\quad n =1,2,\dots\ ,
\ee
where $\alpha^{(n)}$ is an arbitrary element
of the braid group ${\cal B}_n$.
\begin{defin} 
The space ${\cal C}(R,F)$ is a commutative subalgebra in ${\cal M}(R,F)$. It is called the  
characteristic subalgebra of ${\cal M}(R,F)$. 
\end{defin}

The most famous examples of the QM-algebras are given by compatible Yang--Baxter 
matrix pairs $\{R,P\}$ and $\{R,R\}$.

 \medskip
The QM algebras ${\cal M}(R,P)$, related to the series of quantum Lie groups, were introduced in 
\cite{FRT}.
They can be naturally interpreted as noncommutative deformations of algebras of functions on Lie 
groups. 

 \medskip
The QM-algebras 
${\cal M}(R,R)$ are known under the name {\em reflection equation algebras}. Their investigations
were initiated in \cite{Ch,KS}. 
These algebras have special properties among the quantum matrix algebras.
In particular, for the reflection equation algebras, their characteristic subalgebras are central: 
${\cal C}(R,R)\subset Z[{\cal M}(R,R)]$. Geometrically, these algebras can be interpreted as 
deformed algebras of (left or right) adjoint-invariant differential operators on Lie groups.

 \medskip
Three sets of generators for the characteristic subalgebra ${\cal C}(R,F)$ of the $O(k)$ type QM-
algebra
are described in \cite{OP-BMW}:
\begin{enumerate}
\item The set $\{g,\goth{p}_i\}_{i\geq 0}$, where 
\be
\lb{power sums}
\goth{p}_0:=\tr_{\!\! R}\, {\rm Id}_{_V}={ (q-\mu)(q^{-1}+\mu)\over (q-q^{-1}) }, \quad 
\goth{p}_1 := \tr_{\!\! R}\, M\, ,\quad
\goth{p}_i := {\rm  ch}(\sigma_{i-1}\dots\sigma_2\sigma_1) 
,  \; \forall\,i\geq 2.
\ee
generates ${\cal C}(R,F)$ under conditions (\ref{init-cond}) on the parmeters $q$ and $\mu$. For 
reasons that will become clear in section \ref{sec6} elements $\goth{p}_i$ are called {\it power 
sums}.
\item The set $\{g,\goth{h}_i\}_{i\geq 0}$\,: 
\be
\lb{complete symm-f}
\goth{h}_0 :=1, \quad
\goth{h}_i :={\rm  ch}(s^{(i)}), \; \forall\, i\geq 1,
\ee
is defined under stronger restrictions	(\ref{init-cond}) and (\ref{mu-s}) on the parameters.
Its elements $\goth{h}_i$ are called {\it complete  sums}.
\item The set $\{g,\goth{e}_i\}_{i\geq 0}$\,: 
\be
\lb{elementary symm-f}
\goth{e}_0 :=1, \quad
\goth{e}_i :={\rm  ch}(a^{(i)}), \; \forall\, i\geq 1,
\ee
is defined under conditions (\ref{init-cond}) and (\ref{mu-a}).
Its elements  $\goth{e}_i$ are called {\it elementary sums}. Due to finiteness height condition 
(\ref{height-k}) this set is finite: $\goth{e}_{k+i}=0\; \forall\, i>0$. Moreover, 
the rank-one condition
(\ref{rank=1}) leads to algebraic relations among the generators of this set. These dependencies  are 
described in the theorem below.
\end{enumerate}

\begin{theor}\lb{theorem4.9}{\rm\bf \cite{OP-reciprocal}}
Let the QM-algebra ${\cal M}(R,F)$ be of the orthogonal $O(k)$ type. Then
the following reciprocal relations 

\be
\lb{reciprocal} 
g^{k-i}\, {\goth e}_i\, =\, {\goth e}_k\, {\goth e}_{k-i}\,\qquad \forall\; i=0,1,\dots ,k\, 
\ee
are satisfied.
\end{theor}

\begin{rem}\lb{remark1.4} 
{\rm 		
We note, without going into details, that due to the accidental isomorphisms \cite{JO} of the quantum
groups in low dimensions, the case $k=4$  
of the orthogonal type Yang--Baxter matrices is exceptional. Here the word `dimension' means the dimension of the space 
$V$ such that $R^\circ\in {\rm Aut}(V\otimes V)$ (see eq.(\ref{R-st})). 
We say `dimension', not `height', because for the standard deformation the dimension of the space
$V$ coincides with the height of $R^\circ$.

\vskip .2cm
The Yang--Baxter representation
of the antisymmetrizer 
$\rho_{R^\circ} (a^{(2)})$ 
splits into two orthogonal projectors
(the $q$-analogues of the self-dual and anti-self-dual parts), 
\be\lb{qsdaqsd}\rho_{R^\circ} (a^{(2)})=A^+ +A^{-}\  ,\ \ \ \ \ 
A^+ A^{-}=A^{-}  A^+=0\ .\ee
The operators $A^+$ and $A^-$ satisfy relations\,\footnote{For $k\neq4$ the operator 
$\rho_{R^\circ} (a^{(2)})$ does not admit a splitting into two orthogonal projectors satisfying
(\ref{svproj}).}
\be\lb{svproj}
A^{\ \varepsilon}_{1}\, R^{\circ}_{2}\, R^{\circ}_{1}=
R^{\circ}_{2}\, R^{\circ}_{1}\,A^{\ \varepsilon}_{2}\ ,\ \ \ \ \   
A^{\ \varepsilon}_{2}\, R^{\circ}_{1}\, R^{\circ}_{2}=
R^{\circ}_{1}\, R^{\circ}_{2}\,\ A^{\ \varepsilon}_{1}\ ,\ \varepsilon=\pm \ .
\ee
Thus, in this case, the resolution of unity ${\rm Id}_{V\otimes V}$ 
includes four mutually orthogonal projectors instead of
three,
\be \lb{4proj}{\rm Id}_{V\otimes V}=\rho_{R^\circ} (s^{(2)})+A^+ +A^{-}+
\rho_{R^\circ} (c^{(2)})\ ,\ee
and the QM-algebra is defined by a bigger set of relations:
\[ \Pi_\alpha M_{\overline{1}}M_{\overline{2}}=M_{\overline{1}}M_{\overline{2}}\Pi_\alpha\]
for all 4 projectors $\Pi_\alpha$ entering the decomposition of unity (\ref{4proj}). So in this case 
we are dealing with the quotient algebra of the QM-algebra considered above in this subsection
but this does not affect the Cayley--Hamilton identity.

\vskip .2cm
Moreover, in addition to the standard Yang--Baxter matrix $R^\circ$ 
one can construct another Yang--Baxter matrix $R^{\circ\circ}$ with a cubic minimal
characteristic polynomial; its projectors are 
\[ A^{+}\ ,\ A^{-}\ \ \mbox{and}\ \  \rho_{R^{\circ}}(s^{(2)})+\rho_{R^{\circ}}(c^{(2)})\ .\]
The geometry of the corresponding four-dimensional quantum Euclidean and Minkowski spaces for the compatible pair $\{R^{\circ},P\}$  
is considered in \cite{OSWZ,OSWZ2}. 

\vskip .2cm
In the article \cite{JO} the accidental isomorphisms are established only for the FRT-algebra,
that is, the QM-algebra for the compatible matrix pair $\{R^\circ,P\}$.
In the times of writing the article \cite{JO} 
the more general QM-algebras were not yet quite at hand. We suspect that the accidental 
 isomorphisms hold in much more general circumstances, at least, for the compatible pairs 
 $\{R^\circ,R^\circ\}$ as well, but certainly, one should be able to say something for a general compatible pair
 $\{R,F\}$ where $R$ is an arbitrary orthogonal $O(4)$ type Yang--Baxter matrix. We believe that 
 this question deserves further investigation.

\vskip .2cm
We also do not know about any investigation of knot invariants or quantum spin 
chain models which make use of this
exceptional Yang--Baxter matrix $R^{\circ\circ}$ --- this should certainly be an interesting task. 
	}\hfill$\blacktriangleright$
\end{rem}

\subsection{Matrix $\star\,$-product and matrix inversion}\label{subsec4.2}

To define proper analogs of the matrix powers for the quantum matrix $M$ we introduce a  space of 
${\cal M}(R,F)$-valued matrices. 
\begin{defin}
\lb{P-space}
Let ${\cal P}(R,F)$ be a linear subspace of ${\rm End}(V)\otimes {\cal M}(R,F)$
spanned by ${\cal C}(R,F)$-multiples of the identity matrix: 
${\rm Id_{_V}}\, {\zeta} \;\;\, \forall\, {\zeta} \in {\cal C}(R,F)$,
and by matrices of the form
\be
\lb{pow}
M^1 :=\! M, \;\; (M^{\alpha^{(n)}})_{1} :=\! \Tr{2,\dots ,n}\Bigl( M_{\overline 1}
\dots M_{\overline n}\,\rho_R(\alpha^{(n)})\Bigr),\ n =2,3,\dots,
\ee
for any $\alpha^{(n)}\in{\cal B}_n$. 
\end{defin}
The space ${\cal P}(R,F)$
carries a natural structure of a right ${\cal C}(R,F)$--module:
\be
\lb{r-module}
M^{\alpha^{(n)}} {\rm ch} (\beta^{(i)}) =M^{(\alpha^{(n)}\beta^{(i)\uparrow n})}\, \ 
\forall\, \alpha^{(n)}\in {\cal B}_n, \; \beta^{(i)}\in {\cal B}_i\, ,\ n,i=1,2,\dots\, ,
\ee
Here, in the right hand side, the
exponent $\alpha^{(n)}\beta^{(i)\uparrow n}$ belongs to the braid group ${\cal B}_{n+i}$.
Its left factor $\alpha^{(n)}$ is
the image of the element $\alpha^{(n)}\in {\cal B}_n$
under the monomorphism
${\cal B}_{n}\hookrightarrow {\cal B}_{n+i}: \sigma_j\mapsto \sigma_{j}$, and the right factor 
$\beta^{(i)\uparrow n}$ is the image of the  element $\beta^{(i)}\in {\cal B}_i$ under the  monomorphism 
${\cal B}_{i}\hookrightarrow {\cal B}_{n+i}: \sigma_j\mapsto \sigma_{j+n-1}$.
In fact, the formula 
(\ref{r-module}) describes a component-wise multiplication of the matrix 
$M^{\alpha^{(n)}}$ by the element
${\rm ch} (\beta^{(i)})$ from the right.
\begin{defin}
\lb{*product}
We call {\it $\star$-product}
the binary operation 
${\cal P}(R,F) \otimes {\cal P}(R,F)\rightarrow \hspace{-4mm}^\star\hspace{3mm} {\cal P}(R,F)$ 
defined by 
\ba
\nn
({\rm Id_{_V}}\, {\rm ch} (\beta^{(i)})) \star  M^{\alpha^{(n)}} :=
M^{\alpha^{(n)}}  {\rm ch} (\beta^{(i)})&=:&
M^{\alpha^{(n)}}\! \star  ({\rm Id_{_V}}\, {\rm ch} (\beta^{(i)}))  ,
\\[2pt]
\nn
({\rm Id_{_V}}\, {\rm ch} (\alpha^{(n)}))\star ({\rm Id_{_V}}\, {\rm ch} (\beta^{(i)})) &:=
& {\rm Id_{_V}}\,( {\rm ch} (\alpha^{(n)}) {\rm ch} (\beta^{(i)})) ,
\\[2pt]
\lb{MaMb}
M^{\alpha^{(n)}}\! \star  M^{\beta^{(i)}} 
&:=&  M^{(\alpha^{(n)}\star \beta^{(i)})} ,
\ea
\be
\lb{umnB}
\mbox{where}\quad\;\;
\alpha^{(n)}\star \beta^{(i)} :=
\alpha^{(n)}\beta^{(i)\uparrow n} (\sigma_n\dots \sigma_2 \sigma_1\sigma_2^{-1}\dots \sigma_n^{-1})
\in{\cal B}_{n+i} .
\ee
\end{defin}

The $\star$-product defines a structure of an associative $\Bbb C$-algebra on ${\cal P}(R,F)$ (see 
\cite{OP-BMW}, proposition 4.9).\footnote{Note that the formula (\ref{umnB}) defines the associative 
composition of the braids. The verification reduces to checking the equality 
$\tau_k \tau_{k+\ell}=\tau_{\ell}^{\uparrow k}\tau_k$, $k,l>0$ in the braid group 
${\cal B}_{k+\ell+1}$. Here 
$\tau_j:=\sigma_j\dots\sigma_1\sigma_2^{-1}\dots \sigma_j^{-1}$, $j>0$. This verification 
is a not difficult exercise.
} 
We call ${\cal P}(R,F)$ {\em quantum matrix (QM-) powers algebra}.

\medskip
Specifying to the $\star$-product of the quantum matrix $M$ by the arbitrary matrix  
$N\in{\cal P}(R,F)$ we have 
\be
\lb{M*}
M \star  N = M\cdot \phi(N)  \quad \forall N\in {\cal P}(R,F),
\ee
where $\cdot$ denotes the usual matrix multiplication and the linear map $\phi$ is given by 
\be
\lb{phi} 
\phi(N)_1 := \Tr{2} \left( N_{\overline 2} R_{1}\right),
\ee
It will be important for us that the map $\phi$ is invertible. Explicit expression for the inverse map is 
given by \cite{OP-BMW}, Lemma 3.13:
\ba
\lb{phi-inv}
\phi^{-1}(N)_1 &=& \mu^{-2}\,\TR{2}{R_F } \left(  F_{1}^{-1}
N_1 F_{1} (R_F )^{-1}_{1} \right),
\ea
where $R_F $ is the twisted Yang--Baxter matrix (\ref{R_F }).

\begin{rem}
{\rm For the family of the reflection equation algebras ${\cal M}(R,R)$ the map $\phi$ is identical 
and so, the left $\star$-multiplication by the matrix $M$ 
coincides with the usual matrix multiplication. \hfill$\blacktriangleright$
}
\end{rem}

Now we introduce an analogue of the quantum matrix {\em $n$-th power} $M^{\overline n}$ $\,
\forall\,n\geq 0\,$:  
\be
\lb{M^k}
M^{\overline{0}} :={\rm Id}_{_V}\, , \qquad M^{\overline{n}}\, :=\,
\underbrace{M\star  M\star \dots \star  M}_{\mbox{\small $n$ times}}\, =\, 
M^{(\sigma_{n-1}\dots\sigma_2\sigma_{1})}\, \equiv\, M^{(\sigma_1\sigma_2\dots \sigma_{n-1})}.
\ee

For a usual matrix algebra ${\cal M}(P,P)$ the corresponding algebra  ${\cal P}(P,P)$, as the 
${\cal C}(P,P)$-module, is spanned by the usual matrix powers  $M^{\overline n}\equiv M^n$, 
$n\geq 0$. It is easy to check that for a compatible pair $\{R,F\}$,
where $R$ is Hecke type,
the ${\cal C}(R,F)$-module ${\cal P}(R,F)$ is spanned by the quantum 
matrix powers $M^{\overline n}$
as well, see \cite{GPS2}, proposition 2.9.

\vskip .2cm
In the BMW case and, in particular, for the orthogonal QM-algebras the situation is different: 

\begin{prop}
\lb{prop4.11}{\em (\cite{OP-BMW}, proposition 4.11)} 
QM-powers algebra ${\cal P}(R,F)$ of the BMW type QM-algebra ${\cal M}(R,F)$ is commutative 
with respect to the $\star$-product operation.
As the ${\cal C}(R,F)$--module, it is spanned by matrices
$$
M^{\overline{n}}\,  \quad \mbox{and}\quad
\Mt(M^{\overline{n+2}})\, , \quad n=0,1,\dots\, 
$$
Here a ${\cal C}(R,F)$--module map~
$\Mt : {\cal P}(R,F)$ $\rightarrow$ $ {\cal P}(R,F)$ is given by
\be
\lb{Mt}
\Mt (N) := M\cdot \xi(N) \;\;  \forall\, N\in {\cal P}(R,F),
\qquad\mbox{where}
\qquad
\xi(N)_1 := \Tr{2} \left( 
N_{\overline 2}  K_{1}\right).
\ee
\end{prop}

In the definition of the orthogonal $O(k)$ type QM-algebra ${\cal M}(R,F)$ we did not demand the 
invertibility of the contaction $g$. The invertibility is needed if one aims:\vspace{-2.3mm}
\begin{itemize}
	\item[  i)] to resolve the  reciprocal relations;\vspace{-2.3mm}
	\item[ ii)] to introduce the inverse powers of the quantum matrix $M$.\vspace{-2mm}
\end{itemize}  

\vskip .2cm
Below we discuss an extension of the QM-algebras by the inverse contraction $g^{-1}$. We will 
formulate all the results in the most general situation, i.e., in the BMW case. For the subfamily of the 
orthogonal $O(k)$ type QM-algebras, these results  stay valid without any modifications.

\medskip
Let ${\cal M^{^\bullet\!}}(R,F)$ and ${\cal C^{^\bullet\!}}(R,F)$ denote, respectively, localizations of the 
BMW type QM-algebra ${\cal M}(R,F)$ and of its characteristic subalgebra ${\cal C}(R,F)$
by adding the inverse  $g^{-1}$ of the contraction:
$
g^{-1}\, g\, =\, g\, g^{-1}\, =\, 1\, 
$.
The algebras ${\cal M^{^\bullet\!}}(R,F)$ and ${\cal C^{^\bullet\!}}(R,F)$ are further called the {\em  
extended} QM-algebra and its characteristic subalgebra, respectively.

\begin{prop}
\lb{proposition4.11} 
{\em(\cite{OP-BMW}, section 4.5)} The contraction and its inverse satisfy following	permutation relations with the QM-algebra generators
\be
\lb{g-perm}
M\, g\, =\, g\, ( G^{-1}M G)\, , \qquad g^{-1}\, M\, =\, (G^{-1}MG)\, g^{-1}\, ,
\ee	
where the mutually inverse numeric matrices $G^{\pm 1}\in {\rm Aut}(V)$ are given by formulas
\be
\lb{G}
G_1 \, :=\, \tr_{(23)} K_2 F_1^{-1} F_{2}^{-1}\, , \qquad  G_1^{-1}\, =\, \tr_{(23)} F_2 F_1 K_2\, .
\ee
The extended QM-algebra ${\cal M^{^\bullet\!}}(R,F)$ contains
the inverse to the quantum matrix $M$,
\be
\lb{M-inv} 
M^{-1}\, :=\, \mu\, \xi(M)\, g^{-1},\qquad
M\cdot M^{-1}\,=\, {\rm Id}_{_V}\, =\, M^{-1}\cdot M.
\ee 
\end{prop}

\vskip .4cm
The matrix $M^{-1}$ is the inversion of $M$ with respect to the usual matrix product. Inversion with respect 
to the $\star$-product looks differently 
\be
\lb{M-star-inv}
M^{\overline{-1}} :=
\phi^{-1}(M^{-1}),\qquad M^{\overline{-1}}\star M \, =\,{\rm Id}_{_V}\,  =\, M\star M^{\overline{-1}}.
\ee
Consider  extension ${\cal P^{^\bullet\!}}(R,F)$ of the QM-powers  algebra ${\cal P}(R,F)$ of the  BMW type QM-algebra ${\cal M}(R,F)$
by means of repeated $\star$-multiplication by $M^{\overline{-1}}$:
\be
\lb{Minv*N}
M^{\overline{-1}}\star  N \,:=\, \phi^{-1}(M^{-1}\cdot N)\, =:\, N \star  M^{\overline{-1}}
\qquad
\forall\; N\in{\cal P^{^\bullet\!}}(R,F)\, .
\ee
\begin{prop} {\em(\cite{OP-BMW}, section 4.5)}
The extended QM-powers algebra  ${\cal P^{^\bullet\!}}(R,F)$ is associative and commutative with respect to the $\star$-product. It also carries a structure of the right ${\cal C^{^\bullet\!}}(R,F)$-module.
The algebra ${\cal P^{^\bullet\!}}(R,F)$ contains arbitrary integer $\star$-powers of the quantum matrix $M$:
\be
\lb{integer-powers}
M^{\overline{i}}\star M^{\overline{n}}\, =\, M^{\overline{i+n}}\;\; \forall\  i,n\in {\Bbb Z}.
\ee
As the ${\cal C^{\bullet}}(R,F)$--module,   ${\cal P^{\bullet}}(R,F)$ is spanned by matrices
$M^{\overline n}$ and $\Mt (M^{\overline n})$, $n\in {\Bbb Z}$.
\end{prop}

\subsection{Resolution of reciprocal relations}\lb{components}

In this subsection we describe resolution of the reciprocal relations for the extended QM-algebra 
${\cal M^{^\bullet\!}}(R,F)$ of the orthogonal type $O(k)$ (for details and proofs 
see\cite{OP-reciprocal}). 
From now on we use the shorthand notation $O_k(R,F)$ for these algebras, 
${\cal C}_k(R,F)$ for their 
corresponding extended characteristic subalgebras, and ${\cal P}_k(R,F)$ for their extended 
QM-powers algebras.

 \medskip
For the subseries of the orthogonal QM-algebras with even value of their index $k=2\ell$ we assume 
additionally that\,\footnote{Sufficient conditions on the Yang--Baxter matrix pair $(R,F)$ which 
guarantee the validity of the assumption are described in the 
Proposition 4.20 in \cite{OP-reciprocal}. }
\be
\lb{assume-2ell}
\mbox{\em
element	$g^{-\ell}{\goth e}_{2\ell}$ belongs to the center of $O_{2\ell}(R,F)$.
}
\ee
In this situation we introduce two quotient algebras
\be
\lb{quotients-O}
O^{\pm}_{2\ell}(R,F) := O_{2\ell}(R,F)/ \langle g^{-\ell}{\goth e}_{2\ell}\mp 1 \rangle,
\ee
which are called {\em positive} and {\em negative components} of the algebra $O_{2\ell}(R,F)$, respectively. Their
corresponding characteristic subalgebras and QM-powers algebras are defined in the standard way
\be
\lb{quotients-CP}
{\cal C}^{\pm}_{2\ell}(R,F) := {\cal C}_{2\ell}(R,F)/ \langle g^{-\ell}{\goth e}_{2\ell}\mp 1 \rangle,\qquad
{\cal P}^{\pm}_{2\ell}(R,F) := {\cal P}_{2\ell}(R,F)/ \langle{\rm Id_{_V}} (g^{-\ell}{\goth e}_{2\ell}\mp 1) \rangle.
\ee

\begin{theor}
\lb{theorem4.23} (\cite{OP-reciprocal})
In the assumption  (\ref{assume-2ell}) the reciprocal relations (\ref{reciprocal}) in the component algebras
$O^{\pm}_{2\ell}(R,F)$ can be resolved with respect to ${\goth e}_j$, $j\geq \ell$:
\be
\lb{reciprocal2} 
{\goth e}_{\ell +i} = \pm g^i\, {\goth e}_{\ell -i}\,\qquad \forall\; i=0,1,\dots ,\ell\, ,
\ee
respectively. This implies, in particular, that the element ${\goth e}_{2\ell}=\pm g^\ell$ is invertible, and that in the negative component $O^-_{2\ell}$ the element ${\goth e}_\ell$ vanishes: 
${\goth e}_{\ell}=0$.
\smallskip

For the algebras $O_{2\ell -1}(R,F)$ 
the reciprocal relations
(\ref{reciprocal}) have a solution
\be
\lb{reciprocal3} 
{\goth e}_{\ell+i}\, =\, g^{i+1/2} {\goth e}_{\ell-1-i}\,\qquad
\forall\; i=0,1,\dots \ell-1\, ,
\ee
where we denote 
\be
\lb{root-g}
g^{1/2} :=  g^{1-\ell} {\goth e}_{2\ell -1}\quad \Rightarrow\quad g=(g^{1/2})^2,\;\; {\goth e}_{2\ell -1}=(g^{1/2})^{2\ell-1}.
\ee
The latter relations imply that the elements $g^{1/2}$ and ${\goth e}_{2\ell -1}$ are invertible.
\end{theor}

\begin{rem}
{\rm For the orthogonal $O(k)$ type QM-algebras one can introduce notion of a {\em quantum determinant} of the quantum matrix $M$: ${ {\rm det}_q} M := q^{k(k-1)}\, {\goth e}_k$ (see \cite{OP-reciprocal}, section 4.3). As we see, the invertibility of the contraction $g$ implies the invertibility of the quantum determinant ${ {\rm det}_q} M$. }\hfill$\blacktriangleright$	
\end{rem}

\begin{cor}\lb{cor4.24}
In the assumption (\ref{assume-2ell}) characteristic subalgebras 
${\cal C}^\pm_{2\ell}(R,F)$
are generated, respectively, by the sets of elements $g\cup\{{\goth e}_i\}_{i=0}^\ell$ and  $g\cup\{{\goth e}_i\}_{i=0}^{\ell-1}$.

\vskip .2cm
Characteristic subalgebra ${\cal C}_{2\ell -1}(R,F)$
is generated by the set  $g^{1/2}\cup\{{\goth e}_i\}_{i=0}^{\ell-1}$.
\end{cor}

\section{Cayley--Hamilton theorems}\lb{sec5}

In this section we obtain our main result ---
Cayley--Hamilton theorems \ref{theorem5.6} and \ref{theorem5.8} for the extended orthogonal QM-algebras ${O}_k(R,F)$. Preparatory results for these theorems are derived in the subsection \ref{subsec5.1}. These results are obtained in a weaker setting of the BMW type QM-algebras.

\medskip
For readers convenience we remind the conditions on the initial set of data of the QM-algebras which are used in the derivations in these section:
\ba
\nonumber
&\mbox{\em $\{R,F\}$ is a compatible pair of strict skew invertible Yang--Baxter matrices.} 
\\[-2pt]
&\mbox{\em In subsection \ref{subsec5.1} 
 the Yang--Baxter matrix $R$  is of BMW type.}
\lb{cond-RF}&\\[-2pt]
\nonumber
&\mbox{\em In subsection \ref{subsec5.2} and later on  the Yang--Baxter matrix $R$  is 
of $O(k)$ type.}&
\ea

\vskip .1cm
We also remind restrictions on the parameters $q$ and $\mu$ of the BMW type Yang--Baxter matrix $R$: first, these are conditions (\ref{init-cond}) needed for the definition of the BMW algebras ${\cal W}
(q,\mu)$;
second, conditions ensuring the existence of the $q$-antisymmetrizers $a^{(i)}\in {\cal W}_n
(q,\mu)$,~ $1\leq i\leq n$, and, hence, of the wedge powers of matrix $M$:~ 
$i_q M^{a^{(i)}}\!\!\in {\cal P}(R,F)$, see eq.(\ref{La}) below:\,\footnote{Notice, that we escape the unnecessary restriction $q^{2n}\neq 1$ in (\ref{mu-initial})  rescaling $i$-th wedge powers of the quantum matrix $M$ by factors $i_q$.} 
\be
\lb{mu-initial}
q\notin \{0,\pm 1\}, \;\; \mu\notin \{0,\pm q^{\pm 1}\},\;\; \mbox{and}\;\; q^{2i}\neq 1,\;\; \mu\neq -q^{1- 2i}\;\; \forall \, i=2,3,\dots ,n-1.
\ee
For the orthogonal $O(k)$ type Yang--Baxter matrices we impose additional restrictions on the parameter $\mu$.
Namely, starting from subsection \ref{subsec5.2} we specify value of $\mu$:
\be
\lb{mu-secondary}
 \mu=q^{1-k},\;\; \mu\neq -q^{-1-2k}.
\ee
Here the left equality is present in the definition \ref{def-O(k)} of the type $O(k)$ Yang--Baxter matrix, while the right inequality guarantees that $R$ 
is not of the $Sp(2k)$ type (see paragraph above the definition \ref{def-O(k)}).
Condition (\ref{mu-secondary}) in combination with (\ref{mu-initial}) for $n=k+1$ gives
\be
\lb{mu-final}
q\neq 0,\;\; q^{2i}\neq 1\;\;\forall\, i=1,\dots,k,\;\; q^{k+2}\neq -1,\;\; \mu=q^{1-k},
\ee
which will be our set of restrictions on the parameters $q$ and $\mu$ in the formulations of the 
Cayley--Hamilton theorems for the QM-algebras ${O}_k(R,F)$.

\subsection{Basic identities}\lb{subsec5.1}

Consider the set of {\em wedge powers} of the quantum matrix $M$: $M^{a^{(i)}}\!\!\in {\cal P}(R,F)$,
where $a^{(i)}\in {\cal W}_n(q,\mu)$. We introduce two sequences of their {\em descendant} matrices:
\be
\lb{La} 
A^{(m,i)} :=\, i_q\, M^{\overline{m}}  \star  M^{a^{(i)}} , \quad B^{(m+1,i)} :=\,
i_q\, M^{\overline{m}} \star  \Mt (M^{a^{(i)}})\;\;\; \forall\, i,m:\, 1\leq i\leq n,\; m\geq 0  .
\ee
It is  suitable to extend the sequences to the boundary values of their indices  $i=0$:
\ba
\lb{LT0}
A^{(m,0)}\ :=0\ \ \ {\mathrm{and}}\ \ \ \ B^{(m,0)}\ :=\ 0\,\qquad \forall\; m\geq 0\, ,
\ea
and $m=-1$ and $m=0$:
\ba
\lb{AB-boundary} 
 A^{(-1,i)}\, :=\, i_q\, \phi^{-1}\left(
\Tr{2,3,\dots i} M_{\overline{2}}M_{\overline{3}}\dots M_{\overline{i}}\, \rho_R(a^{(i)})
\right)\ , \quad B^{(0,i)}\ :=\ i_q\, \phi^{-1}\bigl(\xi\bigl(M^{a^{(i)}}\bigr)\bigr)\, .
\ea
In case if the contraction $g$ and, hence, quantum matrix $M$ are invertible the descendants 
$A^{(-1,i)}$ and $B^{(0,i)}$ allow uniform expressions:
$$
A^{(-1,i)} = i_q M^{\overline{-1}} \star  M^{a^{(i)}}, \qquad 
B^{(0,i)} = i_q M^{\overline{-1}} \star  \Mt (M^{a^{(i)}}).
$$

In \cite{OP-BMW} two series of recursive relations among the descendants were derived:
\begin{lem}\lb{lemma5.1}
Under conditions  (\ref{cond-RF}), (\ref{mu-initial}), the  matrices $A^{(m-1,i+1)}$, $B^{(m+1,i+1)}: \;0\leq i\leq n-1, \; m\geq 0,$
satisfy recurrent relations
\ba\lb{rek1} A^{(m-1,i+1)} &=& q^i M^{\overline{m}}\, {\goth e}_i\, -\,
A^{(m,i)}\, -\, { \mu q^{2i-1}(q-q^{-1})\over 1+\mu q^{2i-1}}\ B^{(m,i)}\, ,
\\[2pt]
\lb{rek2}
B^{(m+1,i+1)} &=&\Biggl( \mu^{-1}q^{-i} M^{\overline{m}}\, {\goth e}_i\, +\,
{q-q^{-1}\over 1+\mu q^{2i-1}}\ A^{(m,i)}\, -\,B^{(m,i)}\Biggr) g\,  .\ea\end{lem}

In \cite{OP-SpCH} these relations were used to obtain explicit expressions in terms of the matrix powers $M^{\overline{j}}$ for certain descendant matrices. Namely:
\begin{cor}
\lb{corollary5.2}
Relations (\ref{rek1}) and (\ref{rek2}) can be resolved for $A^{(m,i)}$:
$1\leq i\leq n$ and $m\geq i-2$ 
\be
\lb{cor1a}
A^{(m,i)}\;\, =\;\, (-1)^{i-1}\sum_{j=0}^{i-1} (-q)^j\Biggl\{
M^{\overline{m+i-j}} + {1-q^{-2}\over 1+\mu q^{2i-3}}\sum_{r=1}^{i-j-1}M^{\overline{m+i-j-2r}} (q^2 g)^r
\Biggr\} {\goth e}_j\, ,
\ee
and for $B^{(m,i)}$: $1\leq i\leq n$ and $m\geq i$
\ba
\nonumber
B^{(m,i)}&=& (-1)^{i-1}\sum_{j=0}^{i-1} (-q)^j\Biggl\{
\mu^{-1} q^{-2j} M^{\overline{m-i+j}} g^{i-j}\hspace{56mm}
\\[-6pt]
&&\hspace{45mm}
-\,{q^{-1}(1-q^{-2})\over 1+\mu q^{2i-3}}\sum_{r=1}^{i-j-1}M^{\overline{m+i-j-2r}} (q^2 g)^r
\Biggr\} {\goth e}_j\, .
\lb{cor1b}
\ea 
\end{cor}

The expansions (\ref{cor1a}), (\ref{cor1b}) were used in \cite{OP-SpCH} to derive characteristic identities 
for the symplectic QM-algebras. 
However, in the orthogonal case they are not enough: we need similar expansions for the 
descendant matrices $A^{(m,i)}$ and $B^{(m,i)}$  for $m<i$. 
In view of the proposition \ref{prop4.11}, these expansions 
may contain terms $\Mt(M^{\overline{j}})$ and $M^{\overline{j}}$. Below we prepare linear 
combinations of the matrices  $A^{(m,i)}$ and $B^{(m,i)}$, $m<i$, 
such that their expansions do not contain the unwanted terms $\Mt(M^{\overline{j}})$.  

\medskip
Let us introduce $2\mbox{\sc n}\times 2\mbox{\sc n}$ matrices $U^{(i)}$ and $V^{(i)}$, 
$i=1,2,\dots ,n$:
\be
\lb{UV(i)} 
U^{(i)}\! :=  \!
\left( \begin{array}{cc}
\!\!\! - I, &\displaystyle{ \!\!\!\!\!\!\!\!\! {q-q^{-1}\over 1+\mu q^{2i-3}}\, g I }
\\[12pt]
\displaystyle{\!\!\! - {\mu q^{2i-3}(q-q^{-1})\over 1+\mu q^{2i-3}}\, I,} &\!\!\!\!\!\!\!\! -g I
 \end{array}\!\!\right)\! , \;
 V^{(i)}\! :=   {1+\mu q^{2i-3}\over 1+\mu q^{2i-1}} \!
 \left(\begin{array}{cc}\!\! (1+\mu q^{2i-3})\, g I, &\!\!\! - (q-q^{-1})\, g I\!
 \\[12pt]
\!\!\mu q^{2i-3}(q-q^{-1})\, I, &\!\!\! -(1+\mu q^{2i-3}) I\!
 \end{array}\!\!\right)\! ,
 \ee
 where in addition to (\ref{mu-initial}) we assume $\mu\neq -q^{1-2n}$.
These matrices form almost mutually inverse pairs:  
\be
\lb{UinvV}
U^{(i)} V^{(i)}\ = \ (1+\mu q^{2i-5}) g I.
\ee 
Using matrices $U^{(i)}$ one can arrange recurrencies (\ref{rek1}), (\ref{rek2}) into a 
$\mbox{\sc n}\times 2\mbox{\sc n}$ matrix form:
\be
\lb{Nx2N}
X^{(m,i)}\ :=\ \left( A^{(m-1,i)}, B^{(m+1,i)} \right)\ =\  M^{\overline{m}}\left(q^{i-1} I,\, 
{g\over \mu q^{i-1}} I\right) {\goth e}_{i-1} +
\left( A^{(m,i-1)}, B^{(m,i-1)} \right) U^{(i)} .
\ee
Multiplying both sides of this equality by the $2\mbox{\sc n}\times 2\mbox{\sc n}$ matrix
\be
\lb{Qi}
Q^{(i)}\ :=\ V^{(i)}\star \left( 
\begin{array}{cc}
 I, &  0
\\[1pt]
 0, & M^{\overline{2}}
\end{array}
\right)
\ee 
we obtain recursive relation for  the matrices $X^{(m,i)}$:
\be
\lb{X-rek}
X^{(m,i)}\star Q^{(i)} = -q^{-1}(1+\mu q^{2i-3})M^{\overline{m}}\star \left(q^{i-2} I,\, {g\over \mu q^{i-2}} 
M^{\overline{2}}\right) g\, {\goth e}_{i-1}\, +\,
(1+\mu q^{2i-5}) X^{(m+1,i-1)} g.
\ee
Note the important difference between the relation (\ref{X-rek}) and the formulas from the Lemma 
\ref{lemma5.1}: the recurrence on index $i$ in it goes in 
the opposite direction $i\rightarrow i-1$ (assuming that $g$ is invertible). Solving explicitely this 
recurrence one  obtains

\begin{prop}
	\lb{prop5.2}
Under assumptions  (\ref{cond-RF}), parameter restrictions (\ref{mu-initial})  and 
$\mu\neq- q^{1-2n}$ the descendant matrices $X^{(m,i)}$ (\ref{Nx2N})
satisfy relations
\ba
\nn
X^{(m,i)}\star Q^{(i)}\star Q^{(i-1)}\star \dots \star Q^{(i-s+1)}\, = \,\prod_{p=2}^s (1+\mu q^{2(i-p)-1})\,
\Bigl[\, -q^{-1}(1+\mu q^{2(i-s)-1})\, M^{\overline{m}}\Bigr. 
\\[2pt]
\lb{X-rel-gen}
\textstyle\star\,\Bigl( q^{i-s-1}{\textstyle \sum_{j=0}^{s-1}}\,(-q)^j M^{\overline{s-j-1}}\, {\goth e}_{i+j-s}\, g^s,
\,\, {1\over \mu q^{i-s-1}}\sum_{j =0}^{s-1}\, (-{1\over q})^{j} M^{\overline{s+j+1}}\, {\goth e}_{i+j-s}\, g^{s-j}
\Bigr)\;
\\[4pt]
\nn
-\,\, (q-q^{-1})\,q^{i-s-2}M^{\overline{m}}\star \sum_{j=1}^{s-1}\sum_{r=j}^{s-1}(-q)^{2r-j-s+1} 
M^{\overline{2(s-r)+j-1}}\, {\goth e}_{i-j}\,g^r 
\Bigl(I,\,\, -q g I\Bigr)\qquad\;
\\[4pt]
\nn
\Bigl. +\,\, (1+\mu q^{2(i-s)-3})\, X^{(m+s,i-s)} g^s\, \Bigr],\qquad\quad\;\;\;\,\forall\; m\geq 0,\;1\leq 
i\leq n,\; 1\leq s\leq i.\qquad\;\;\;
\ea
\end{prop}

\noindent {\bf Proof.}~ First, let us comment on the parameter restrictions in the statement of the 
proposition.
Conditions (\ref{cond-RF}), (\ref{mu-initial}) guarantee correctness of expressions for $X^{(m,i)}$, 
$1\leq i\leq n$, and for $Q^{(j)}$, $1\leq j\leq n-1$.
The additional condition $\mu\neq- q^{1-2n}$ is necessary to avoid singularities in $Q^{(n)}$.
\medskip

We will prove the formula (\ref{X-rel-gen}) by induction on $s$. 

 \smallskip
For $s=1$, it reproduces equality (\ref{X-rek}), thereby setting the induction base.
Carrying out the inductive step $s\rightarrow s+1$ is a straightforward but time-consuming 
calculation. 
We will briefly comment on its main stages.

\vskip .1cm
The right hand side of eq.(\ref{X-rel-gen}) is a sum of three components: the first one occupies two 
upper lines in formula (\ref{X-rel-gen}), the second and third components are written, 
correspondingly, in the third and fourth lines. 
Let us denote these first, second and last summands by symbols 
$\alpha_s$, $\beta_s$ and $\gamma_s$, 
respectively.
Our notation shows explicitly a value of the index $s$ in (\ref{X-rel-gen}).
To prove the inductive step we have to verify the equality
\be\lb{abc}
(\alpha_s + \beta_s +\gamma_s)\star Q^{(i-s)} = \alpha_{s+1} +\beta_{s+1}+ \gamma_{s+1}.
\ee
We carry this calculation separately for the summands.

\medskip
For $\gamma_s$ we have
$$
\gamma_s\star Q^{(i-s)} = \gamma_{s+1} + \alpha_{s+1}|_{_{j=0}}.
$$
Here symbol $\alpha_{s+1}|_{_{j=0}}$ denotes the component with index $j=0$ in the sum 
$\sum_{j=0}^s$ in the expression for $\alpha_{s+1}$ (see the second line in (\ref{X-rel-gen})).

\medskip
In the calculation with $\beta_s$ we use the formula
$$
\Bigl(I, -q g I\Bigr)\star Q^{(i-s)} = (1+\mu q^{2(i-s)-3}) \Bigl(-I,\, q^{-1}M^{\overline{2}}\Bigr)\, g
$$
and obtain
$$
\beta_s\star Q^{(i-s)}\, =\, \Bigl( \beta^{(1)}_{s+1}|_{_{r>j}},\, \beta^{(2)}_{s+1}|_{_{r<s}}\Bigr).
$$
Here $\beta^{(1)}_{s+1}$ and $\beta^{(2)}_{s+1}$ denote the left and right 
$\mbox{\sc n}\times \mbox{\sc n}$ matrix components of the $\mbox{\sc n}\times 2\mbox{\sc n}$ 
matrix $\beta_{s+1}$;
notation $|_{_{r>j}}\,$,$\; |_{_{r<s}}$ describe additional restrictions on the summation index $r$ in 
$\beta^{(1)}_{s+1}$ and $\;\beta^{(2)}_{s+1}$.

\medskip
For the calculation with $\alpha_{s}$ it is suitable to separate 
the matrix $Q^{(i-s)}$ into block-diagonal and 
block-antidiagonal parts $Q^{(i-s)}=Q_d^{(i-s)}+Q_a^{(i-s)}$:
$$
Q^{(i-s)}_d:=\,{\textstyle -\, {(1+\mu q^{2(i-s)-3})^2\over 1+\mu q^{2(i-s)-1}} }\!\left( 
\!\!\begin{array}{cc}
g I, &\!\!\!\!  0
\\[1pt]
0, &\!\!\!\! M^{\overline{2}}
\end{array}\!\!
\right)\!,\quad
Q^{(i-s)}_a:=\, {\textstyle {(q-q^{-1})(1+\mu q^{2(i-s)-3})\over 1+\mu q^{2(i-s)-1}} }\!\left( \!\!
\begin{array}{cc}
0 , &\!\!\!\!  -g M^{\overline{2}}
\\[1pt]
\mu q^{2(i-s)-3} I, &\!\!\!\! 0
\end{array}\!\!
\right)\!.
$$
We obtain
$$
\alpha_{s}\star Q_d^{(i-s)} =  \alpha_{s+1}|_{_{j\neq 0}}, \qquad \alpha_{s}\star Q_a^{(i-s)}=  
\Bigl( \beta^{(1)}_{s+1}|_{_{j=r}},\, \beta^{(2)}_{s+1}|_{_{r=s}}\Bigr).
$$
Collecting the results of all calculations we finally prove the induction assertion (\ref{abc}).
\hfill$\blacksquare$\medskip

Later in this section  we derive the Cayley--Hamilton identities for the QM algebras in the orthogonal 
case.
It is common for three series -- the linear \cite{IOP1}, the symplectic \cite{OP-SpCH} and the 
orthogonal QM-algebras   -- that these identities originate from 
the finiteness height property of their defining Yang--Baxter 
matrix $R$. However, details of the derivation are 
different for each particular series.
\medskip
 
For the BMW type Yang--Baxter matrix $R$, assuming it has height $k$ (see eq.(\ref{height-k})) one 
immediately obtains two sequences of basic matrix identities in the algebra
${\cal P}(R,F)$
\be
\lb{ahah} 
A^{(m-1,k+1)}\equiv 0\, , \quad B^{(m,k+1)}\equiv 0\, \quad \forall\; m\geq 0\, .
\ee
In general, the left hand sides of these identities can be expanded into 
linear combinations of the terms 
$M^{\overline{j}}$ and $\Mt(M^{\overline{j}})$,
$j=0,1,\dots$ . Our next goal is to prepare a combination of the identities (\ref{ahah}) which contains 
only terms $M^{\overline{j}}$ in the expansion.

\medskip
For the symplectic QM-algebras of the type $Sp(2k)$ the  right equalities in (\ref{ahah}) follow from 
the left ones (see \cite{OP-SpCH}, remark 4.8). 
Here one obtains the Cayley--Hamilton identity by the application of the corollary 
\ref{corollary5.2} to the relation  $A^{(k-1,k+1)}\equiv 0$.   
The resulting characteristic identity is polynomial of degree $2k$. 

\medskip
For the orthogonal type $O(k)$ the situation is different: the left and right identities in eq.(\ref{ahah}) 
are independent and one can take their combinations
to obtain characteristic identities of the lower degree. 
We have prepared such combinations in 
the proposition \ref{prop5.2}. They are
\be
\lb{ahah1}
X^{(0,k+1)}\star Q^{(k+1)}\star Q^{(k)}\star \dots \star Q^{(k-s+2)}|_{s=\lfloor {k+1\over 2}\rfloor}\, = 
\, 0.
\ee
It is left to calculate explicit expressions for the matrix polynomials in the l.h.s. of (\ref{ahah1}). The 
calculation is different for the
cases of even and odd $k$  and we describe them separately in the next two subsections.

\subsection{Type $O(k)$ Cayley--Hamilton theorem for even  $k= 2\ell$}\lb{subsec5.2}

From now on we consider orthogonal QM-algebras ${O}_k(R,F)$
imposing the conditions (\ref{cond-RF}) on the Yang--Baxter matrices $R$ and $F$, 
and the restrictions (\ref{mu-final})
on the parameters $q$ and $\mu$.
These restrictions guarantee that the matrix expression standing in the left hand side of eq.
(\ref{ahah1}) is correctly defined and
that one can use proposition \ref{prop5.2} for its calculation.
\medskip

First, we consider the case of even $k$: $k=2\ell$, $\ell=1,2,\dots\; .$ 
\smallskip

Noticing that the common factor $\prod_{p=2}^s (1+\mu q^{2(i-p)-1})$ in the 
r.h.s. of (\ref{X-rel-gen}) is invertible 
we can rewrite the identity  (\ref{ahah1}) in the following form:
\ba
\nn
0\; =\; -\, (q+q^{-1})\,
\Bigl( q^{\ell} \sum_{j=0}^{\ell-1}\,(-q)^j M^{\overline{\ell-j-1}}\, {\goth e}_{\ell+j+1}\, g^\ell,\,\, q^{\ell-1}\sum_{j 
=0}^{\ell-1}\, (- q)^{-j} M^{\overline{\ell+j+1}}\, {\goth e}_{\ell+j+1}\, g^{\ell-j}\Bigr)
\\[2pt]
\lb{Xeven}
-(q-q^{-1})\, q^{\ell -1} \sum_{j=1}^{\ell -1}\sum_{r=j}^{\ell -1}(-q)^{2r-j-\ell +1} M^{\overline{2(\ell -r)+j-1}}\, {\goth e}_{2\ell -j+1}\,g^r 
\Bigl(I,\,\, -q g I\Bigr)\; +\; 2 X^{(\ell,\ell +1)} g^\ell\, .
\ea
Next, we can express the matrix components of $X^{(\ell,\ell +1)}= \Bigl( A^{(\ell -1, \ell +1)},\, B^{(\ell +1,\ell +1)}\Bigr)$ via matrix powers of $M$ with the help of 
corollary \ref{corollary5.2}. Namely, substituting $\mu=q^{1-2\ell}$ into eqs.(\ref{cor1a}), (\ref{cor1b}) 
we find
\ba
\lb{Aeven}
A^{(\ell -1, \ell +1)} &\!\!=&\!\! (-1)^{\ell}\sum_{j=0}^{\ell} (-q)^j\Bigl\{
M^{\overline{2\ell -j}} + {1-q^{-2}\over 2}\sum_{r=1}^{\ell -j}M^{\overline{2(\ell-r) -j}} (q^2 g)^r
\Bigr\} {\goth e}_j ,
\\[2pt]
\lb{Beven}
B^{(\ell +1,\ell +1)} &\!\!=&\!\! (-1)^{\ell }\sum_{j=0}^{\ell} (-q)^j\Bigl\{ q^{2\ell -2j-1} M^{\overline{j}} g^{\ell +1-j}\,
	-\,{1-q^{-2}\over 2 q}\sum_{r=1}^{\ell -j}M^{\overline{2(\ell-r) -j+2}} (q^2 g)^r
\Bigr\} {\goth e}_j .
\ea   

Consider the left and right matrix components of the identity (\ref{Xeven}) separately. 
Upon substitution of the expression (\ref{Aeven}) for $A^{(\ell -1, \ell +1)}$ the left matrix component 
assumes the form
\ba
\nn
0& =& (\lambda - 2 q)\,
 q^{\ell} \sum_{j=0}^{\ell-1}\,(-q)^j M^{\overline{\ell-j-1}}\, {\goth e}_{\ell+j+1}\, g^\ell
 \\
 \nn
 && -\;\lambda\, q^{\ell -1} \sum_{j=1}^{\ell -1}\sum_{r=j}^{\ell -1}(-q)^{2r-j-\ell +1} M^{\overline{2(\ell -r)+j-1}}\, {\goth e}_{2\ell -j+1}\,g^r 
 \\
 \lb{AA3}
 &&+\; (-1)^{\ell}\sum_{j=0}^{\ell} (-q)^j\Bigl\{2\,
 M^{\overline{2\ell -j}}\, +\, {\lambda\over q}\,\sum_{r=1}^{\ell -j}M^{\overline{2(\ell-r) -j}} (q^2 g)^r
 \Bigr\} {\goth e}_j g^\ell .
\ea
Here we introduced a notation 
\be
\lb{lambda}
\lambda := q-q^{-1}
\ee
to simplify further transformations of the expression. 

\medskip
Let us combine the term with the coefficient $\lambda$ from the first 
line of eq.(\ref{AA3}) with the part of the
term with $\lambda$ from the third line that corresponds to the value $\ell -j$ of the summation index 
$r$. Making the 
substitution $j\leftrightarrow \ell -j-1$ of the summation index in the first of these terms we obtain
\be
\lb{ex1}
\lambda (-1)^\ell\, \sum_{j=0}^{\ell-1}\,
(-q)^{2\ell-j-1} M^{\overline{j}}\Bigl( {\goth e}_{2\ell -j} g^j- {\goth e}_j g^\ell\Bigr) g^{\ell -j} .
\ee
Combining the rest of terms with the coefficient $\lambda$ from the third line of eq.(\ref{AA3}) with 
the term from the second line we obtain
\be
\lb{ex2}
\lambda (-1)^\ell\, \sum_{j=0}^{\ell-2}\, \sum_{r=1}^{\ell-j-1}
(-q)^{2r+j-1} M^{\overline{2(\ell-r)-j}}\Bigl( {\goth e}_{2\ell -j} g^j- {\goth e}_j g^\ell\Bigr) g^{r} .
\ee
Obviously, expressions (\ref{ex1}) and (\ref{ex2}) can be joined together
\be
\lb{ex12}
(-1)^{\ell+1}(1-q^{-2})\, \sum_{j=0}^{\ell-1}\, \sum_{r=1}^{\ell-j}
(-q)^{2r+j} M^{\overline{2(\ell-r)-j}}\Bigl( {\goth e}_{2\ell -j} g^j- {\goth e}_j g^\ell\Bigr) g^{r}\, =:\, (-1)^{\ell+1} U_\ell \, .
\ee
Notation $U_\ell$ is introduced here for the future convenience.

\medskip
Finally, we observe that the rest of terms from the first and the last lines in eq.(\ref{AA3}) can be unified into a single sum
\be
\lb{ex3}
2 (-1)^\ell\, \sum_{j=0}^{2\ell}\, (-q)^j M^{\overline{2\ell -j}} {\goth e}_j g^\ell .
\ee 
Combining expressions (\ref{ex12}) and (\ref{ex3}) and rescaling the result we obtain a simpler expression for the identity
(\ref{AA3})
\be
\lb{exAfin}
\sum_{j=0}^{2\ell}\, (-q)^j M^{\overline{2\ell -j}} {\goth e}_j g^\ell
- {1\over 2}\, U_\ell
\, =\, 0\, .
\ee

Now let us consider the right matrix component in the identity (\ref{Xeven}). Substituting the 
expression 
(\ref{Beven}) for $B^{(\ell +1,\ell +1)}$ into it we get
\ba
\nn
0& =& -\, (\lambda +  2 q^{-1})\,
q^{\ell-1} \sum_{j=0}^{\ell-1}\,(-q)^{-j} M^{\overline{\ell+j+1}}\, {\goth e}_{\ell+j+1}\, g^{\ell -j}
\\
\nn
&& +\;\lambda\, q^{\ell} \sum_{j=1}^{\ell -1}\sum_{r=j}^{\ell -1}(-q)^{2r-j-\ell +1} M^{\overline{2(\ell -r)+j-1}}\, {\goth e}_{2\ell -j+1}\,g^{r+1} 
\\
\lb{BB3}
&&+\, (-1)^{\ell }\sum_{j=0}^{\ell} (-q)^j\Bigl\{ 2\, q^{2\ell -2j-1} M^{\overline{j}} g^{\ell +1-j}\,
-\,{\lambda\over  q^2}\sum_{r=1}^{\ell -j}M^{\overline{2(\ell-r) -j+2}} (q^2 g)^r
\Bigr\} {\goth e}_j g^\ell .
\ea
Here again, we used the notation $\lambda$ (\ref{lambda}) to bring the identity to a form suitable for  
transformations.

\medskip
Unifying the term with $\lambda$ from the first line of eq.(\ref{BB3}) with the part of the 
term with $\lambda$ 
from the third line that corresponds to $r=1$ 
we obtain
\be
\lb{exB1}
\lambda (-1)^\ell\, \sum_{j=0}^{\ell-1}\,
(-q)^j M^{\overline{2\ell -j}}\Bigl( {\goth e}_{2\ell -j} g^j- {\goth e}_j g^\ell\Bigr) g.
\ee
Combining the rest of terms with $\lambda$ from the third line of eq.(\ref{BB3}) with the term from 
the second line we obtain
\be
\lb{exB2}
\lambda (-1)^\ell\, \sum_{j=0}^{\ell-2}\, \sum_{r=2}^{\ell-j}
(-q)^{2r+j-2} M^{\overline{2(\ell-r)-j+2}}\Bigl( {\goth e}_{2\ell -j} g^j- {\goth e}_j g^\ell\Bigr) g^{r} .
\ee
Expressions (\ref{exB1}) and (\ref{exB2}) can be joined together
\be
\lb{exB12}
q^{-1} (-1)^\ell\,(1-q^{-2}) \sum_{j=0}^{\ell-1}\, \sum_{r=1}^{\ell-j}
(-q)^{2r+j} M^{\overline{2(\ell-r)-j+2}}\Bigl( {\goth e}_{2\ell -j} g^j- {\goth e}_j g^\ell\Bigr) g^{r}\, =\, q^{-1} (-1)^\ell\, 
M^{\overline{2}} \star U_\ell\, .
\ee
Finally, the rest of terms from the first and the last lines in eq.(\ref{BB3}) can be unified into a single 
sum
\be
\lb{exB3}
2\, q^{-1}(-1)^\ell\, \sum_{j=0}^{2\ell}\, (-q)^{2\ell - j} M^{\overline{j}} {\goth e}_j g^{2\ell-j+1} .
\ee 
Combining expressions (\ref{exB12}) and (\ref{exB3}) and rescaling the result we obtain a concise 
version of the identity
(\ref{BB3}) 
\be
\lb{exBfin}
\sum_{j=0}^{2\ell}\, (-q)^j M^{\overline{2\ell -j}} {\goth e}_{2\ell -j} g^{j+1}
+ {1\over 2}\, M^{\overline{2}} \star U_\ell
\, =\, 0\, .
\ee

For further simplification of identities (\ref{exAfin}) and (\ref{exBfin}) we derive alternative 
expressions for $U_\ell$ and $M^{\overline{2}}\!\star U_\ell$.
Opening the bracket $(1-q^{-2})$ in  (\ref{ex12}), we represent the matrix $U_\ell$ as a difference of 
two double sums.
Changing the summation index $r\rightarrow r-1$ in the second double sum (the one with the 
common
factor $-q^{-2}$), one can bring the expression for $U_\ell$ to the form
\ba
\nonumber
U_\ell \; =\; - \,{\displaystyle \sum_{j=0}^{2\ell} (-q)^j M^{\overline{2\ell-j}}}
( {\goth e}_{2 \ell -j} g^j - {\goth e}_j g^{\ell})\,
\hspace{64mm} 
\\
\lb{altU}
\hspace{25mm} 
+\, (M^{\overline{2}}- g I)\star
{\displaystyle \sum_{j=0}^{\ell -1}\, \sum_{r=0}^{\ell-j-1} } (-q)^{j+2r}
M^{\overline{2(\ell -r)-j-2}} ({\goth e}_{2\ell-j}g^j - {\goth e}_j g^{\ell}) g^r .
\ea
In a similar way we transform the expression for $M^{\overline{2}}\!\star U_\ell$:
\ba
\nonumber 
M^{\overline{2}} \star  U_{\ell} \; =\;  -\sum_{j=0}^{2\ell} (-q)^j M^{\overline{2\ell-j}}
({\goth e}_{2\ell-j}g^j - {\goth e}_j g^{\ell}) \, g 
\hspace{70mm}
\\
\lb{altMU}
\hspace{30mm} 
+\, q^2 (M^{\overline{2}}- g I)\star
\sum_{j=0}^{\ell -1}\, \sum_{r=0}^{\ell-j-1} (-q)^{j+2r}
M^{\overline{2(\ell-r)-j-2}} ({\goth e}_{2\ell-j}g^j - {\goth e}_j g^{\ell}) \, g^{r+1} .
\ea
Substituting the expressions (\ref{altU}) and (\ref{altMU}), respectively, into eqs. (\ref{exAfin}) and 
(\ref{exBfin}) 
and comparing the results we conclude that, under the  
assumptions that contraction $g$ is invertible and 
$q^2\neq -1$ (the latter is guaranteed by conditions (\ref{mu-final})), the union of
eqs. (\ref{exAfin}) and 
(\ref{exBfin}) is equivalent to the 
following two matrix equalities
\ba
\lb{ortho-5} 
&&\sum_{j=0}^{2\ell}(-q)^j M^{\overline{2\ell-j}} ({\goth e}_{2\ell -j} g^j + {\goth e}_j g^{\ell}) \; =\;  0\, ,
\\
\lb{ortho-6} 
&&
\hspace{-22mm}
(M^{\overline{2}}-gI)\star  \sum_{j=0}^{2\ell -2} (-q)^j\, M^{\overline{2\ell-j-2}}
\hspace{-7mm}\sum_{r={\rm max}\{0,j+1-\ell\}}^{\lfloor j/2\rfloor}\hspace{-8mm}
({\goth e}_{2\ell-j+2r}g^{j-2r} - {\goth e}_{j-2r}g^{\ell})\, g^r\; =\;  0\, .
\ea
Here the double sum in the last formula is obtained from the double sum in eq.(\ref{altMU}) by a 
substitution of the summation index $j \rightarrow j+2r$.
\medskip

Up to now, we didn't apply the
rank-one property 
specific for the orthogonal Yang--Baxter matrices. 
This property gives rise to the reciprocal relations among the elements ${\goth e}_j$ (\ref{reciprocal}) which 
allows one
to factorize identities (\ref{ortho-5}) and (\ref{ortho-6}) 
\ba
\lb{ortho-7} 
&&\Bigl\{\sum_{j=0}^{2\ell} (-q)^j M^{\overline{2\ell -j}} {\goth e}_j\Bigr\} ({\goth e}_{2\ell}+g^{\ell})\, =\, 0\, ,
\\ 
\lb{ortho-8}
&&\Bigl\{(M^{\overline{2}}-gI)\star \sum_{j=0}^{2\ell-2} (-q)^j\, M^{\overline{2\ell-j-2}}
\hspace{-7mm}\sum_{r={\rm max}\{0,j+1-\ell\}}^{\lfloor j/2\rfloor} \hspace{-8mm}
{\goth e}_{j-2r}g^r\Bigr\} ({\goth e}_{2\ell}-g^{\ell})\, =\,  0\, .
\ea
Relations (\ref{ortho-7}) and (\ref{ortho-8}) are $q$-analogues of the characteristic identities which are 
fulfilled by the quantum matrix $M$ of the orthogonal QM-algebras of the type $O(2\ell)$. They 
assume more familiar form in the extended QM-algebras $O_{2\ell}(R,F)$. Indeed, by theorem
\ref{theorem4.23}, in the assumption (\ref{assume-2ell}), 
the left factors in formulas (\ref{ortho-7}) and 
(\ref{ortho-8}) represent the characteristic identities in the components $O^{+}_{2\ell}(R,F)$ and 
$O^{-}_{2\ell}(R,F)$, respectively.
Thus, we arrive to the following $q$-analogue of the Cayley--Hamilton theorem:

\begin{theor}\lb{theorem5.6}
Consider  the extended orthogonal QM-algebra $O_{2\ell}(R,F)$.
Assume that the parameters $q$ and $\mu$ of the orthogonal type Yang--Baxter 
matrix $R$ satisfy conditions 
(\ref{mu-final}) for $k=2\ell$. Assume additionally
that the element $g^{-\ell}{\goth e}_{2\ell}\in O_{2\ell}(R,F)$ is central, and hence, the two quotient
algebras of $O_{2\ell}(R,F)$ -- the components  $O^{\pm}_{2\ell}(R,F)$ (see (\ref{quotients-O})) -- 
are well defined.

\medskip
Then the quantum matrix $M$ of generators of the positive  component $O^+_{2\ell}(R,F)$ satisfies  
the identity
\be
\lb{CH-O+} 
\sum_{j=0}^{2\ell}(-q)^j\, M^{\overline{2\ell -j}}\, {\goth e}_j\, =\, 0,
\ee
where ${\goth e}_{\ell+j}= g^j {\goth e}_{\ell -j}$, $\forall\; j=1,2,\dots ,\ell$.

\medskip
The quantum matrix $M$ of 
generators of the negative component $O^-_{2\ell}(R,F)$ satisfies the identity
\be
\lb{CH-O-1}
\left(M^{\overline{2}} - g I\right)\star
\sum_{j=0}^{2\ell-2}(-q)^j\, M^{\overline{2\ell -j-2}}\,\epsilon_j\, =\, 0,
\ee
where
\be
\lb{CH-O-2}
\epsilon_j\, :=\, \sum_{r=0}^{\lfloor j/2\rfloor } {\goth e}_{j-2r}\, g^r\, , \qquad
\epsilon_{\ell -1+j}\, :=\, \epsilon_{\ell-1-j}\, g^{j}\,\qquad \forall\; j=0,1,\dots ,\ell -1 .
\ee
\end{theor}
\noindent {\bf Proof.}~ The Cayley--Hamilton  identities (\ref{CH-O+}) and (\ref{CH-O-1}) for the 
components $O^+(2\ell)$ and $O^-(2\ell)$ follow 
immediately from eqs.(\ref{ortho-7}) and (\ref{ortho-8}). 
\hfill$\blacksquare$

\subsection{Type $O(k)$ Cayley--Hamilton theorem for odd $k=2\ell -1$}\lb{subsec5.3}

In this subsection we continue to investigate the identity (\ref{ahah1}) for the QM-algebras 
$O_{k}(R,F)$, 
now odd $k=2\ell -1$, $\ell=2,\dots$ . We impose the
conditions (\ref{cond-RF}) on the Yang--Baxter matrices $R$ and $F$, and 
the restrictions (\ref{mu-final})
on the parameters $q$ and $\mu$, so that we can use the results 
of proposition \ref{prop5.2}.

\smallskip
First, we rewrite the identity (\ref{ahah1}) using the formula (\ref{X-rel-gen}) and removing the
 nonzero common factor in it:
\ba
\nn
0\; =\; 
\Bigl( q^{\ell-1} \sum_{j=0}^{\ell-1}\,(-q)^j M^{\overline{\ell-j-1}}\, {\goth e}_{\ell+j}\, g^\ell,\,\, q^{\ell-1}\sum_{j =0}^{\ell-1}\, (- q)^{-j} M^{\overline{\ell+j+1}}\, {\goth e}_{\ell+j}\, g^{\ell-j}\Bigr)\hspace{15mm}
\\[2pt]
\lb{Xodd}
+\,(1-q^{-1})\, q^{\ell -1} \sum_{j=1}^{\ell -1}\sum_{r=j}^{\ell -1}(-q)^{2r-j-\ell +1} M^{\overline{2(\ell -r)+j-1}}\, {\goth e}_{2\ell -j}\,g^r 
\Bigl(I,\,\, -q g I\Bigr)\; -\;   X^{(\ell,\ell)} g^\ell\, .
\ea
Then, we express the matrix components of $X^{(\ell,\ell)}= \Bigl( A^{(\ell -1, \ell)},\, B^{(\ell +1,\ell)}\Bigr)$ via matrix powers of $M$ using formulae (\ref{cor1a}), (\ref{cor1b}):
\ba
\lb{Aodd}
A^{(\ell -1, \ell)} &\!\!=&\!\! (-1)^{\ell-1}\sum_{j=0}^{\ell-1} (-q)^j\Bigl\{
M^{\overline{2\ell -j-1}} + (1-q^{-1})\sum_{r=1}^{\ell -j-1}M^{\overline{2(\ell-r) -j-1}} (q^2 g)^r
\Bigr\} {\goth e}_j ,
\\[2pt]
\lb{Bodd}
B^{(\ell +1,\ell)} &\!\!=&\!\! (-1)^{\ell-1 }\sum_{j=0}^{\ell-1} (-q)^j\Bigl\{ q^{2\ell -2j-2} M^{\overline{j+1}} g^{\ell -j}\,
-\,{1-q^{-1}\over  q}\sum_{r=1}^{\ell -j-1}M^{\overline{2(\ell-r) -j+1}} (q^2 g)^r
\Bigr\} {\goth e}_j .
\ea   
Substituting expression (\ref{Aodd}) for $A^{(\ell -1, \ell)}$ into the left matrix component of (\ref{Xodd}) we obtain
\ba
\nn
0& =& 
q^{\ell-1} \sum_{j=0}^{\ell-1}\,(-q)^j M^{\overline{\ell-j-1}}\, {\goth e}_{\ell+j}\, g^\ell
\\
\nn
&& +\;( \underline{\underline{ q^{\ell -1}}} - q^{\ell-2}) \sum_{j=1}^{\ell -1}\sum_{r=j}^{\ell -1}(-q)^{2r-j-\ell +1} 
M^{\overline{2(\ell -r)+j-1}}\, {\goth e}_{2\ell -j}\,g^r 
\\
\lb{CC3}
&&+\; (-1)^{\ell}\sum_{j=0}^{\ell-1} (-q)^j\Bigl\{
M^{\overline{2\ell -j-1}}\, +\, (1-\underline{\underline{q^{-1}}})\,\sum_{r=1}^{\ell -j-1}M^{\overline{2(\ell-r) -j-1}} (q^2 g)^r
\Bigr\} {\goth e}_j g^\ell .
\ea
Combining the terms  with doubly underlined factors from the second and third lines of (\ref{CC3}) we obtain
\be
\lb{Cex1}
 (-1)^{\ell-1}\, \sum_{j=1}^{\ell-1}\, \sum_{r=0}^{\ell-j-1}
(-q)^{j} M^{\overline{2(\ell-r)-j-2}} (q^2g)^r\star\Bigl( M {\goth e}_{2\ell -j} g^{j-1}- {\goth e}_{j-1} g^\ell I\Bigr)g .
\ee
Unifying rest of terms in the third line we get
\be
\lb{Cex2}
 (-1)^{\ell}\, \sum_{j=0}^{\ell-1}\, \sum_{r=0}^{\ell-j-1}
(-q)^{j} M^{\overline{2(\ell-r)-j-1}} (q^2g)^r {\goth e}_j g^\ell  .
\ee
Combining the rest of the terms in the second line with the term from the first line we obtain
\be
\lb{Cex3}
(-1)^{\ell-1}\, \sum_{j=0}^{\ell-1}\, \sum_{r=0}^{\ell-j-1}
(-q)^{j} M^{\overline{2(\ell-r)-j-2}} (q^2g)^r {\goth e}_{2\ell-j-1} g^{j+1}  .
\ee
Expressions (\ref{Cex2}) and (\ref{Cex3}) can be joined together:
\be
\lb{Cex4}
(-1)^{\ell-1}\, \sum_{j=0}^{\ell-1}\, \sum_{r=0}^{\ell-j-1}
(-q)^{j} M^{\overline{2(\ell-r)-j-2}} (q^2g)^r  \star\Bigl({\goth e}_{2\ell-j-1} g^{j} I- M {\goth e}_{j} g^{\ell-1}\Bigr) g .
\ee
Finally, unifying expressions (\ref{Cex1}) and (\ref{Cex4}) and rescaling the result we rewrite the 
identity (\ref{CC3}) in a simplified form:
\be
\lb{CC4} 
\sum_{j=0}^{\ell -1}\sum_{r=0}^{\ell-j-1} (-q)^{j} M^{\overline{2(\ell-r)-j-2}}(q^2g)^r\star\Bigl(
({\goth e}_{2\ell-j-1} g^j I - M {\goth e}_j g^{\ell-1} ) +  (M {\goth e}_{2\ell -j} g^{j-1}  - {\goth e}_{j-1}g^{\ell} I)\Bigr)g\,=\, 0 .
\ee
Here we have extended the range of summation over 
$j$ in the expression (\ref{Cex1}) assuming ${\goth e}_{-1}={\goth e}_{2\ell}:=0$.

\medskip
Now we turn to a simplification of the right matrix component of the identity (\ref{Xodd}). Substituting expression (\ref{Bodd}) for $B^{(\ell +1,\ell )}$ into it we get
\ba
\nn
0\!& =&\! 
q^{\ell-1} \sum_{j=0}^{\ell-1}\,(-q)^{-j} M^{\overline{\ell+j+1}}\, {\goth e}_{\ell+j}\, g^{\ell-j}
\\
\nn
\!&&\! +\,( q^{\ell -1} -  \underline{\underline{q^{\ell}}}) \sum_{j=1}^{\ell -1}\sum_{r=j}^{\ell -1}(-q)^{2r-j-\ell +1} 
M^{\overline{2(\ell -r)+j-1}}\, {\goth e}_{2\ell -j}\,g^{r+1}
\\
\lb{DD3}
\!&&\!+\, (-1)^{\ell}\sum_{j=0}^{\ell-1} (-q)^j\Bigl\{
\underline{q^{2\ell -2j-2}M^{\overline{j+1}}g^{\ell-j}}\, +\, (\underline{\underline{q^{-2}}}-q^{-1})\,\sum_{r=1}^{\ell -j-1}M^{\overline{2(\ell-r) -j+1}} (q^2 g)^r
\Bigr\} {\goth e}_j g^\ell .\qquad
\ea
Combining the  term from the first line of (\ref{DD3}) with the part of the underlined term in 
the third line 
corresponding to the value $j=\ell-1$ of the summation index, 
and the part of the term with the doubly underlined factor $q^{-2}$ in the 
same line  corresponding to the value $r=1$ of the summation index, we obtain
\be
\lb{Dex3}
(-1)^{\ell}\, \sum_{j=0}^{\ell-1}
(-q)^{j} M^{\overline{2\ell-j-1}} \star \Bigl({\goth e}_{j} g^{\ell} I- M {\goth e}_{2\ell-j-1} g^{j}\Bigr) g .
\ee
Combining the term with the doubly underlined factor from the second line of (\ref{DD3}) 
with the rest of the underlined 
and doubly underlined terms from the third line we get
\be
\lb{Dex4}
(-1)^{\ell}\, \sum_{j=0}^{\ell-2}\,\sum_{r=1}^{\ell-j-1} (-q)^{j} 
M^{\overline{2(\ell-r)-j-1}} (q^2 g)^r \star \Bigl({\goth e}_{j} g^{\ell} I- M {\goth e}_{2\ell-j-1} g^{j}\Bigr) g .
\ee
Expressions (\ref{Dex3}) and (\ref{Dex4}) can be joined together:
\be
\lb{Dex5}
(-1)^{\ell}\, \sum_{j=0}^{\ell-1}\,\sum_{r=0}^{\ell-j-1} (-q)^{j} 
M^{\overline{2(\ell-r)-j-1}} (q^2 g)^r \star \Bigl({\goth e}_{j} g^{\ell} I- 
M {\goth e}_{2\ell-j-1} g^{j}\Bigr) g .
\ee
The rest of terms from the second and third lines in eq.(\ref{DD3}) can be unified into 
\be
\lb{Dex6}
(-1)^{\ell}\, \sum_{j=1}^{\ell-1}\,\sum_{r=0}^{\ell-j-1} (-q)^{j} 
M^{\overline{2(\ell-r)-j-1}} (q^2 g)^r \star \Bigl(M {\goth e}_{j-1} g^{\ell} - {\goth e}_{2\ell-j} g^{j} I\Bigr) g .
\ee
Unifying expressions (\ref{Dex5}) and (\ref{Dex6}) and rescaling the result we rewrite the identity (\ref{DD3}) in the form:
\be
\lb{DD4} 
\sum_{j=0}^{\ell -1}\sum_{r=0}^{\ell-j-1} (-q)^{j} M^{\overline{2(\ell-r)-j-1}}(q^2g)^r\star\Bigl(
({\goth e}_{j} g^\ell I - M {\goth e}_{2\ell-j-1} g^{j} ) +  (M {\goth e}_{j-1} g^{\ell}  - {\goth e}_{2\ell-j}g^{j} I)\Bigr)g\,=\, 0 .
\ee

Now we take into consideration reciprocal relations (\ref{reciprocal})  specific for the type 
\mbox{$O(2\ell-1)$} orthogonal QM-algebras.
First of all we observe that the two matrix identities (\ref{CC4}) and (\ref{DD4}) are not independent. 
Their left hand sides satisfy equality
$$
M\star \Bigl( \mbox{l.h.s. of (\ref{CC4})}\Bigr) {\goth e}_{2\ell-1}\, =\, \Bigl( \mbox{l.h.s. of (\ref{DD4})}\Bigr) g^{\ell -1} .
$$
In the extended QM-algebra ${O}_{2\ell-1}(R,F)$ (i.e., if the contraction $g$ is invertible) these two identities are equivalent and can be written in the form 
\be
\lb{ortho-9}
\sum_{j=0}^{\ell -1}\sum_{r=0}^{\ell-j-1}  (-q)^{j} M^{\overline{2(\ell-r)-j-2}}
(q^2g)^r\star \Bigl( ({\goth e}_{2\ell-1} I - M g^{\ell-1} )  {\goth e}_j + (M {\goth e}_{2\ell-1}
 - g^{\ell} I){\goth e}_{j-1}\Bigr)\, =\, 0\, .
\ee
Thus, we obtain the following $q$-analogue of the Cayley--Hamilton theorem

\begin{theor}\lb{theorem5.8}
Consider  the extended orthogonal QM-algebra $O_{2\ell-1}(R,F)$.
Assume that the parameters $q$ and $\mu$ of the orthogonal type Yang--Baxter matrix $R$ satisfy conditions 
(\ref{mu-final}) for $k=2\ell-1$.
Then the quantum matrix $M$ of generators of  $O_{2\ell-1}(R,F)$ satisfies the identity
\be
\lb{CH-O-odd}
\left(M - g^{1/2} I\right)\star \sum_{j=0}^{2\ell-2}(-q)^j\, M^{\overline{2\ell -j-2}}\,
\epsilon_j\, =\, 0\, ,
\ee
where
\be
\lb{CH-O-odd2}
\epsilon_j\, :=\, \sum_{r=0}^{j} {\goth e}_{j-r}\, (-g^{1/2})^r\, , \qquad
\epsilon_{\ell -1 +j}\, :=\, \epsilon_{\ell-1-j}\, g^j\,\qquad\forall\; j=0,1,\dots ,\ell -1\, .
\ee
\end{theor}
\noindent {\bf Proof.}~ In the extended QM-algebra  $O_{2\ell-1}(R,F)$ the reciprocal relations 
admit the  explicit solution (see theorem \ref{theorem4.23}).
The corresponding expressions for elements ${\goth e}_j$ are given in eqs.(\ref{reciprocal3}). Substituting 
these expressions into
(\ref{ortho-9}) and collecting coefficients in the matrix powers $M^{\overline{j}}$ we obtain 
upon cancelling the apperaing common factor 
$(-g^{\ell -1})$ the identity (\ref{CH-O-odd}) with the coefficients given by the formulas 
(\ref{CH-O-odd2}).
\hfill$\blacksquare$

\begin{rem}
{\rm The classical group $O_{2\ell -1}$,  as the classical group $O_{2\ell}$, has two connected 
components. We note that the QM-algebra  $O_{2\ell-1}(R,F)$ has the analogues of the two 
components $O_{2\ell-1}^\pm (R,F)$ as well, which we illustrate in the situation  
when the elements $g$ and ${\goth e}_{2\ell-1}$ are central (see propositions 4.6, 4.12 and 4.21 in 
\cite{OP-reciprocal} for the condititions of the centrality of these elements).
In this case we can fix $g=1$, and, due to the reciprocal relation 
${\goth e}_{2\ell -1}^2=g^{2\ell-1}$, we are left with two possibilities for the choice of 
${\goth e}_{2\ell -1}$, ${\goth e}_{2\ell -1}=\pm 1$. These two choices correspond to the two 
components $O_{2\ell-1}^\pm (R,F)$. The element $g^{1/2}$ (see the formula  (\ref{root-g})) 
in the Cayley--Hamilton identity (\ref{CH-O-odd}) 
takes the values $g^{1/2} :=\pm 1$ for the components $O_{2\ell-1}^\pm (R,F)$, respectively.
 }	\hfill$\blacktriangleright$
\end{rem}

\section{Spectral parameterization }\lb{sec6}

In this section we continue investigation of the extended orthogonal  QM-algebras ${O}_k(R,F)$. We 
assume that the conditions of theorems 
\ref{theorem5.6} and \ref{theorem5.8} are fulfilled. In particular, two 
components $O^\pm_{2\ell}(R,F)$ of the even $k=2\ell$ algebra $O_{2\ell}(R,F)$ are well defined. 
\begin{defin}\lb{E2p}
Let	${\cal E}_{2n}$ denote a $\Bbb C$-algebra
  of polynomials in $2n+1$ pairwise commuting and invertible variables
$\nu_j$, $~j=0,1,\dots ,2n$, satisfying conditions
\be
\lb{specO}
\nu_{n+j}\, \nu_{n+1-j}\, =\, \nu_0^2\,\quad\forall\; j=1,2,\dots ,n\, .
\ee
\end{defin}

Later in this section the variables $\nu_j$ will be given the role of eigenvalues of the orthogonal 
quantum matrix $M$, and therefore we will call them {\em  spectral variables}. 

\medskip
In subsection \ref{subsec6.1} we introduce three homomorphic maps $\pi^\pm_{2\ell}$ and 
$\pi_{2\ell-1}$ from the extended  characteristic subalgebras 
 of the QM-algebras ${O}_{k}(R,F)$ into the algebras of polynomials in the spectral variables.
Using these maps we define central extensions  ${\cal P}^{\pm\, ext}_{2\ell}(R,F)$ and 
${\cal P}^{ext}_{2\ell-1}(R,F)$ of the extended QM-powers algebras ${\cal P}^\pm_{2\ell}(R,F)$ and 
${\cal P}_{2\ell-1}(R,F)$. 
Then the Cayley--Hamilton identities from theorems \ref{theorem5.6} and 
\ref{theorem5.8}, considered as 
relations in the central extensions of the QM-powers algebras,  acquire completely factorized forms. 

\medskip
We stress that, in general, we cannot construct the factorization of the 
Cayley--Hamilton identities in the QM-algebra ${O}_k(R,F)$. However in the situation when the
characteristic subalgebra belongs to the center of the algebra 
${O}_k(R,F)$ we can use the maps $\pi^\pm_{2\ell}$ and 
$\pi_{2\ell-1}$ to define the central extensions of the whole QM-algebra (and not only of the
 QM powers algebra). In this case, the 
Cayley--Hamilton identities take the factorized form in the 
central extension of the QM-algebra. The 
important case of such situation is the RE-algebra ${O}_k(R,R)$.

\medskip
In the final subsection \ref{subsec6.2} we obtain parameterizations for the other two generating sets 
of the characteristic subalgebra -- power sums ${\goth p}_i$ and complete sums 
${\goth h}_i$ - in terms of the spectral variables. 

\medskip
Our results in theorem \ref{corollary5.7} and proposition \ref{corollary6.4} are consistent with the 
results of
\cite{Mudr}, section 8.3,  although they are obtained  by a different method and in a more general 
setting.

\subsection{Factorization of the Cayley--Hamilton identities}\lb{subsec6.1} 

Let us define the mappings of the sets $\{g,\goth{e}_i\}_{i\geq 0}$ (see 
 eq.(\ref{elementary symm-f}))
of the generators of the characteristic subalgebras ${\cal C}^+_{2\ell}(R,F)$, 
${\cal C}^-_{2\ell}(R,F)$ and ${\cal C}_{2\ell-1}(R,F)$ into the polynomial algebras in the spectral 
variables:
\begin{itemize}
\item
 $\pi^+_{2\ell}$ maps the generators of the characteristic subalgebra ${\cal C}^+_{2\ell}(R,F)$ 
 into 
 ${\cal E}_{ 2\ell}$
 \be
 \lb{rep-charO+}
 	\ g\mapsto \nu_0^2 , \;\; {\goth e}_i\mapsto e_i(\nu_1, \nu_2,\dots ,\nu_{2\ell})\,
 	\quad \forall\; i=1,2,\dots 2\ell ;
\ee
\item
 $\pi^-_{2\ell}$ maps the generators of the characteristic subalgebra ${\cal C}^-_{2\ell}(R,F)$ 
 into ${\cal E}_{ 2\ell-2}$
 \be
 \lb{rep-charO-}
 	\ g\mapsto \nu_0^2 , \;\; {\goth e}_i\mapsto e_i(\nu_0,-\nu_0,\nu_1, \nu_2,\dots ,\nu_{2\ell-2})\,
 	\quad \forall\; i=1,2,\dots 2\ell;
 \ee
\item
 $\pi_{2\ell-1}$ maps the generators of  the characteristic subalgebra 
 ${\cal C}_{2\ell-1}(R,F)$  
  into ${\cal E}_{ 2\ell-2}$
 \be\lb{rep-charOodd}
 g^{1/2}\mapsto \nu_0 ,\;\;\; {\goth e}_i\mapsto
 e_i(\nu_0,\nu_1, \nu_2,\dots ,\nu_{2\ell-2})\,\quad \forall\; i=1,2,\dots 2\ell-1 .
 \ee
\end{itemize}
Here $e_i(\dots)$ are the elementary symmetric polynomials in their arguments.

\begin{theor}
\lb{corollary5.7}
In the settings of theorems \ref{theorem5.6} and \ref{theorem5.8} let us additionally assume that
the elements of the generating sets $g\cup\{{\goth e}_i\}_{i=0}^\ell$, 
$g\cup\{{\goth e}_i\}_{i=0}^{\ell-1}$ and 
$g^{1/2}\cup\{{\goth e}_i\}_{i=0}^{\ell-1}$ of the characteristic subalgebras ${\cal C}^+_{2\ell}(R,F)$, ${\cal C}^-_{2\ell}(R,F)$ and ${\cal C}_{2\ell-1}(R,F)$, respectively,  
are algebraically independent.

\medskip	
Then  the maps $\pi^+_{2\ell}$, $\pi^-_{2\ell}$ and $\pi_{2\ell-1}$ extend to the homomorphisms of 
the respective characteristic subalalgebras into polynomial algebras in the spectral variables. These 
homomorphisms endow the algebras ${\cal E}_{2\ell}$ and ${\cal E}_{2\ell-2}$
with the natural (left) module structure over the corresponding  characteristic subalgebras.
Consider the following central extensions of the QM-powers algebras:
\ba
\mbox{for $O^+_{2\ell}(R,F)$:}	 && {\cal P}^{+\, ext}_{2\ell}(R,F)\, :=\, {\cal P}^+_{2\ell}(R,F)\raisebox{-4pt}{$\bigotimes\atop
{\cal C}^+_{2\ell}(R,F)$} {\cal E}_{2\ell};
\\[2pt]
\mbox{for $O^+_{2\ell}(R,F)$:} &&	{\cal P}^{-\, ext}_{2\ell}(R,F)\, :=\,
{\cal P}^-_{2\ell}(R,F)\raisebox{-4pt}{$\bigotimes\atop {\cal C}^-_{2\ell}(R,F)$} {\cal E}_{2\ell-2};
\\[2pt] 
\mbox{for $O_{2\ell-1}(R,F)$:} &&	{\cal P}^{\,ext}_{2\ell -1}(R,F)\, :=\,
{\cal P}_{2\ell -1}(R,F)\raisebox{-4pt}{$\bigotimes\atop
{\cal C}_{2\ell -1}(R,F)$} {\cal E}_{2\ell -2}.
\ea
In these centrally extended QM-powers algebras the Cayley--Hamilton
identities (\ref{CH-O+}), (\ref{CH-O-1}) and (\ref{CH-O-odd}) acquire completely factorized forms, 
respectively,
\ba
\lb{CHO+factor} 
\mbox{for~ ${\cal P}^{+\, ext}_{2\ell}(R,F)$: }	&& 
{\prod_{i=1}^{2\ell}}\!\star \!\left(
M - q\nu_i I \right)\, =\, 0,
\\
\lb{CHO-factor} 
\mbox{for~ ${\cal P}^{-\, ext}_{2\ell}(R,F)$: }
&&
(M-\nu_0 I)\star (M+\nu_0 I){\prod_{i=1}^{2\ell-2}} \!\!\star\! \left( M - q\nu_i I
\right)\, =\, 0,
\\
\lb{CHOodd-factor}
\mbox{for~ ${\cal P}^{\,ext}_{2\ell -1}(R,F)$: }
&& (M-\nu_0 I ){\prod_{i=1}^{2\ell-2}}\!\!
\star\!\left( M - q\nu_i I\right)\, =\, 0 .
\ea
For the subfamily of the extended orthogonal reflection equation algebras  these matrix equalities 
are valid, respectively, in the algebras ${O^+_{2\ell}(R,R)}$, ${O^-_{2\ell}(R,R)}$ and 
${O_{2\ell-1}(R,R)}$ with the replacement of the $\star$-product by the usual matrix product.
\end{theor}

\begin{rem}{\rm
For the classical  orthogonal groups, the functions ${\goth e}_i$, $i=1,\dots ,\ell$,
on the two  connected components of $O_{2\ell}$: $O^+_{2\ell}\equiv SO_{2\ell}$ and 
$O^-_{2\ell+2}$,  and on $O_{2\ell+1}$ are
algebraically independent. This justifies, at least perturbatively, the 
assumption about the independence of the elements ${\goth e}_i$ in the theorem.
}\hfill$\blacktriangleright$
\end{rem}

\noindent {\bf Proof of theorem \ref{corollary5.7}.}
The verification that the mappings (\ref{rep-charO+})-(\ref{rep-charOodd}) allow 
homomorphic extensions
consists in checking the reciprocal relations (\ref{reciprocal2}), (\ref{reciprocal3}) for their images. 
These calculations become almost obvious if one takes into account the following relations for the 
elementary symmetric polynomials
$$
e_{k-i}(\nu_1,\dots ,\nu_k) = e_i(\nu_1^{-1},\dots ,\nu_k^{-1}) \,e_k(\nu_1, \dots ,\nu_k)\quad 
\forall\, i=0,1,\dots, k ,
$$
and then uses the equalities

$$
e_i(\nu_1^{-1},\dots ,\nu_k^{-1}) = \nu_0^{-2i}\, e_i(\nu_1, \dots ,\nu_k),
$$
which are valid in case if the set of the arguments $\{\nu_j\}_{j=1}^k$ is invariant under the involutive map $\nu_j^{-1}\leftrightarrow {\nu_j/{\nu_0^2}}\;\;\; \forall j=1,\dots ,k$.
\smallskip

The factorization of the Cayley--Hamilton identity (\ref{CH-O+}) for the positive component $O^+_{2\ell}(R,F)$ is obvious. 

To verify the factorization of the Cayley--Hamilton identity (\ref{CH-O-1}) for the negative component $O^-_{2\ell}(R,F)$ it suffices to use the following property of the elementary symmetric polynomials.
Denote $\{\nu\}_k:= \{\nu_1,\nu_2,\dots ,\nu_k\}$. Then~
$
e_i(\{\nu\}_k) = e_i(\nu_0,-\nu_0,\{\nu\}_k) + \nu_0^2\, e_{i-2}(\{\nu\}_k),
$~
and hence,
$$
e_i(\{\nu\}_k) = 
\sum_{r=0}^{\lfloor i/2\rfloor}
e_{i-2r}(\nu_0,-\nu_0,\{\nu\}_k)\, \nu_0^{2r}\qquad \mbox{(c.f. with eqs.(\ref{CH-O-2})).}
$$

To verify the factorization of the Cayley--Hamilton identity (\ref{CH-O-odd}) for the algebra $O_{2\ell-1}(R,F)$ one applies the similar property: $e_i(\{\nu\}_k) = e_i(\nu_0,\{\nu\}_k) - \nu_0\, e_{i-1}(\{\nu\}_k),$ and hence,
$$
e_i(\{\nu\}_k) = 
\sum_{r=0}^{i}
e_{i-r}(\nu_0,\{\nu\}_k) (-\nu_0)^{r}\qquad \mbox{(c.f. with eqs.(\ref{CH-O-odd2})).}
$$
\hfill$\blacksquare$

\subsection{Parameterization of Newton and Wronski relations}\label{subsec6.2}

In this final subsection we derive the spectral parameterizations for
two other sets of the generating elements of the characteristic subalgebras in 
the extended orthogonal QM-algebras, namely, for the sets  
$\{{\goth h}_i,g\}_{i=0,1,\dots}$ and $\{{\goth p}_i,g\}_{i=0,1,\dots}$. These elements were introduced 
in subsection \ref{subsec4.1}.

\medskip
We will use the  $(q,\mu)$-analogues of Newton and Wronski relations
obtained for the generators of the characteristic subalgebras of the BMW type QM-algebras in 
\cite{OP-BMW}, theorem 5.2:   
\ba
\lb{Newton-a}
&&\sum_{i=0}^{n-1} (-q)^i {\goth e}_i\, {\goth p}_{n-i} \,=\, (-1)^{n-1} n_q\, {\goth e}_n \, +\, (-1)^n q
\sum_{i=1}^{\lfloor {n/2}\rfloor}\Bigl( q^{n-2i-k} -q^{-n+2i}\Bigr)\, {\goth e}_{n-2i}\, g^i,
\\
\lb{Newton-s}
&&\sum_{i=0}^{n-1} q^{-i} {\goth h}_i\, {\goth p}_{n-i} \,=\,  n_q\, {\goth h}_n \, +
\,\sum_{i=1}^{\lfloor {n/2}\rfloor}\Bigl(
q^{n-2i-1}+q^{-n+2i+1-k}\Bigr)\, {\goth h}_{n-2i}\, g^i,
\\
\lb{Wronski}
&&\sum_{i=0}^n (-1)^i {\goth e}_i\, {\goth h}_{n-i}\, =\, \delta_{n,0} -\delta_{n,2}\, g .
\ea
Here $n$ is an arbitrary positive integer and
$\delta_{i,j}$ is a Kronecker symbol.
We also restrict $\mu=q^{1-k}$ in the formulas of \cite{OP-BMW} since in this text  we are dealing 
with the orthogonal $O(k)$ type QM-algebras.
\begin{rem}{\rm 
Notice that  the elements ${\goth e}_i$ and ${\goth h}_i$, $i=1,\dots ,n$, were originally defined, 
respectively, under the restrictions (\ref{mu-a}) and (\ref{mu-s}) on the parameters $q$ and $\mu$.
However, afterwards we can use the  formulas (\ref{Newton-a}), (\ref{Newton-s}) for an iterative 
construction of these elements in terms of ${\goth p}_i$, starting from ${\goth e}_0={\goth h}_0=1$. 
The latter is possible under weaker restrictions
\be
\lb{weaker-cond}\hspace{5.5cm}
q\neq 0, \quad q^{2i}\neq 1\quad \forall i=2,\dots , n.\hspace{4.5cm}\blacktriangleright
\ee
}\hfill
\end{rem}

Now we derive spectral parameterization for the elements ${\goth h}_i$ and ${\goth p}_i$ using 
relations (\ref{Newton-a})-(\ref{Wronski}) and the parametrization 
formulas (\ref{rep-charO+})-(\ref{rep-charOodd}) for ${\goth e}_i$. 

\begin{theor}
\lb{corollary6.4} In the setting of theorem \ref{corollary5.7} 
assume additionally that the parameter $q$ of the $O(k)$ type QM-algebra satisfies the
conditions (\ref{weaker-cond}) for some integer $n\geq k$.\footnote{For $n=1,2,\dots , k$ the 
inequalities (\ref{weaker-cond}) are fulfilled due to the conditions (\ref{mu-final}) imposed in theorems 
\ref{theorem5.6}, \ref{theorem5.8} and \ref{corollary5.7}.} Adopting the notation from subsections 
\ref{subsec5.2}, \ref{subsec5.3} we redenote $k=2\ell$,  or $k=2\ell-1$ depending on the parity of 
$k$. 

\vskip .2cm
Then under the homomorfisms $\pi^+_{2\ell}$,  $\pi^-_{2\ell}$ and $\pi_{2\ell-1}$ the generating 
elements ${\goth h}_i$ and ${\goth p}_i$, $i=1,2,\dots ,n$, of the characteristic subalgebras	of the 
extended orthogonal type 
QM-algebras, respectively, $O^+_{2\ell}(R,F)$, $O^-_{2\ell}(R,F)$ and $O_{2\ell-1}(R,F)$ have the 
following images
\ba
\lb{para-s-O+} 
\pi^+_{2\ell}:&&\hspace{-4mm} {\goth h}_i \mapsto h_i(\nu_1,\dots ,\nu_{2\ell})
- \nu_0^2\, h_{i-2}(\nu_1,\dots ,\nu_{2\ell}) ,\qquad
{\goth p}_i \mapsto q^{i-1}\sum_{j=1}^{2\ell} d_j (\nu_j)^i ,
\\[-2pt]
\lb{para-p-O+}
&&\hspace{-4mm}
\mbox{where}\;\;
d_j\, :=\!\!\!\! \prod _{r=1\atop r\neq j,\, 2\ell +1-j}^{2\ell}\!\!\!\!\!
(\nu_j - q^{-2}\nu_r)/( \nu_j - \nu_r);
\\[4pt]
\lb{para-s-O-} 
\pi^-_{2\ell}:&&\hspace{-4mm} {\goth h}_i \mapsto h_i(\nu_1,\dots ,\nu_{2\ell-2}) , \qquad
{\goth p}_i \mapsto q^{i-1}\sum_{j=1}^{2\ell-2} d_j (\nu_j)^i + \left(1 + (-1)^i\right) q^{1-2\ell} (\nu_0)^i ,
\\[-2pt]
\lb{para-p-O-}
 &&\hspace{-4mm} \mbox{where}\quad
d_j\, :=
{(\nu_j - q^{-4}\nu_{2\ell -1-j})/ (\nu_j - \nu_{2\ell -1-j})}\!\!\!\!
\prod_{r=1\atop r\neq j,\, 2\ell -1-j}^{2\ell -2}\!\!\!\! {(\nu_j - q^{-2}\nu_r)/ (\nu_j - \nu_r)} ;
\\[2pt]
\lb{para-s-O-odd}
\pi_{2\ell-1}:&&\hspace{-4mm} {\goth h}_i \mapsto h_i(\nu_1,\dots ,\nu_{2\ell-2})
+\nu_0\, h_{i-1}(\nu_1,\dots ,\nu_{2\ell-2}) ,\;\; {\goth p}_i \mapsto q^{i-1}\sum_{j=1}^{2\ell-2} d_j (\nu_j)^i +
q^{2-2\ell}(\nu_0)^i ,\qquad
\\[-2pt]
\lb{para-p-O-odd}
&&\hspace{-4mm}\mbox{where} \quad
d_j\, := {(\nu_j - q^{-2}\nu_0)/(\nu_j - \nu_0)}\!\!\!\!
\prod _{r=1\atop r\neq j,\, 2\ell -1-j}^{2\ell -2}\!\!\!\!
{(\nu_j - q^{-2}\nu_r)/( \nu_j - \nu_r)} .
\ea
Here, in the expressions for ${\goth h}_i$ symbols
$h_i(\dots)$ denote the complete symmetric polynomials in their arguments. In the formulas for 
${\goth p}_i$ we assume that the algebras of polynomials in the spectral variables are extended by 
all the 
inverse differences $(\nu_{r_1}-\nu_{r_2})$ which appear in the quantities $d_j$. 
\end{theor}

\begin{rem}
{\rm Although the expressions for the power sums $\goth{p}_i$, $i=0,1,2,\ldots$,
 with the use of the quantities $ d_j$ from 
theorem \ref{corollary6.4}
require the invertibility of all the differences $(\nu_{r_1}-\nu_{r_2})$, effectively the elements   
$\goth{p}_i$ are polynomials in the variables $q,\, q^{-1}$ and in the spectral variables 
$\nu_0,\nu_1,\nu_2,...$ This follows immediately from the relations (\ref{Newton-a}), understood
as a recursive definition of the power sums.}~\hfill$\blacktriangleright$
\end{rem}

\noindent {\bf Proof of theorem \ref{corollary6.4}.}~
For the proof, we use the following  result.
\begin{lem}
\lb{lemma6.5} Let the elements ${\goth e}_i$, ${\goth h}_i$, ${\goth p}_i$, $i=0,1,\dots ,n,$ satisfy 
the recursive relations 
(\ref{Newton-a})-(\ref{Wronski}) and the initial conditions ${\goth e}_0={\goth h}_0=1$,
${\goth p}_0=q^{1-k}(1+(k-1)_q)$.\footnote{The initial condition for ${\goth p}_0$ 
corresponds to the choice  $\mu=q^{1-k}$ in the definition (\ref{power sums}), i.e., to the specialization to the case of the orthogonal $O(k)$ type QM-algebras.} 
Consider their modifications  ${\goth e}'_i$, ${\goth h}'_i$, 
${\goth p}'_i$ and ${\goth p}''_i$:
\ba
\lb{mod-as}
&&{\goth e}'_i = {\goth e}_i + {\goth e}'_{i-2}\, g , \qquad {\goth h}'_i = {\goth h}_i + {\goth h}'_{i-2}\, g\quad \forall\; i\geq 0 , \quad \mbox{where}
\quad {\goth e}'_{i}={\goth h}'_i=0 \mbox{~~if~~} i<0 ;
\\[2pt] 
\lb{mod-p1}
&&{\goth p}'_0=q^{-k} k_q ,  \qquad\qquad {\goth p}'_1={\goth p}_1 , \qquad\qquad {\goth p}'_i = {\goth p}_i + (q^2\,{\goth p}'_{i-2}-{\goth p}_{i-2})\, 
g\quad\; \forall\; i\geq 2;
\\[2pt]
\lb{mod-p2} 
&&{\goth p}''_0=q^{2-k} (k-2)_q ,  \quad
{\goth p}''_1={\goth p}_1, \qquad\qquad {\goth p}''_i = {\goth p}_i + (q^{-2}{\goth p}''_{i-2}-{\goth p}_{i-2})\, g\;\;\; \forall\; i\geq 2.
\ea
The modified elements ${\goth e}'_i$, ${\goth h}'_i$, ${\goth p}'_i$ and ${\goth p}''_i$, $i=0,1,\dots ,n$, satisfy 
the simpler recursions
\ba
\lb{mod-N} 
&&\sum_{i=0}^{n-1} (-q)^i {\goth e}_i\, {\goth p}'_{n-i}\, =\,  (-1)^{n-1} n_q\, {\goth e}_n\, ,
\qquad \sum_{i=0}^{n-1} q^{-i} {\goth h}_i\, {\goth p}''_{n-i}\, =\,  n_q\, {\goth h}_n, 
\\
\lb{mod-W}
&& \sum_{i=0}^{n} (-1)^i {\goth e}'_i\, {\goth h}_{n-i}\, =\,\delta_{n,0}\, =\,\sum_{i=0}^{n} (-1)^i {\goth e}_i\, {\goth h}'_{n-i}\,. 
\ea	
\end{lem}\noindent {\bf Proof.}~
Let us prove first  the Wronski relations (\ref{mod-W}).
For $n=0,1,$ they are clearly satisfied. For $n\geq 2$, we apply induction on $n$.
E.g., the left equality in (\ref{mod-W}) follows  from the calculation
$$\sum_{i=0}^n (-1)^i {\goth e}'_i\, {\goth h}_{n-i}\, =\,\sum_{i=0}^n (-1)^i ({\goth e}_i + {\goth e}'_{i-2}\, g)\, {\goth h}_{n-i}\, =\,
\delta_{n,0} -\delta_{n,2}\, g +\sum_{i=0}^{n-2} (-1)^i {\goth e}'_i\, {\goth h}_{n-2-i}\, g
\, =\, \delta_{n,0}\, ,$$
where the relations (\ref{mod-as}) and (\ref{Wronski}) and the induction assumption were 
successively used.
The right equality in (\ref{mod-W}) is proved in the same way.

\medskip
To verify  the  relations (\ref{mod-N})
we also apply induction on $n$. Consider, for instance, the left equality in (\ref{mod-N}). Cases 
$n=1$ and $n=2$ providing the induction base are checked directly. Here in the latter case we use 
equality
\be
\lb{popo}
q^2 {\goth p}'_0 - {\goth p}_0\,=\, q - q^{1-k} .
\ee
For $n\geq 3$ we calculate
\ba
\nn
&&\sum_{i=0}^{n-1} (-q)^i {\goth e}_i\, {\goth p}'_{n-i}\, =\, \sum_{i=0}^{n-1} (-q)^i {\goth e}_i\, {\goth p}_{n-i} + \sum_{i=0}^{n-3} (-q)^i {\goth e}_i
\Bigl(q^2 {\goth p}'_{n-i-2} - {\goth p}_{n-i-2}\Bigr) g + (-q)^{n-2} {\goth e}_{n-2}\Bigl(q^2 {\goth p}'_0 - {\goth p}_0 \Bigr)g
\\[2pt]
\nn
&&=
(-1)^n\, q
\sum_{i=1}^{\lfloor {n/2}\rfloor}\Bigl( q^{n-2i-k} -q^{-n+2i}\Bigr) {\goth e}_{n-2i}\, g^i \,
-\,(-1)^{n-2} \, q \!\!\! \sum_{i=1}^{\lfloor {n/2}\rfloor-1}\!\Bigl( q^{n-2-2i-k} -q^{-n+2+2i}\Bigr) {\goth e}_{n-2-2i}\, g^{i+1}
\\
\nn
&&
\;\;\,+\,  
(-1)^{n-1}\, n_q \,{\goth e}_n + (q^2-1)(-1)^{n-3}(n-2)_q\, {\goth e}_{n-2}\,g\,  + \, (-q)^{n-2}(q-q^{1-k})\,{\goth e}_{n-2}\,g
\\[3pt]
\nn
&&=  (-1)^{n-1}\, n_q \,{\goth e}_n\, + \,(-1)^{n}\Bigl\{-(q^2-1)(n-2)_q+(q^{n-1}-q^{n-1-k})+ (q^{n-1-k}-q^{3-n})\Bigr\} {\goth e}_{n-2}\, g
\\[3pt]
\nn
&&=  (-1)^{n-1}\, n_q \,{\goth e}_n.
\ea
Here in the first line we applied formulas (\ref{mod-p1}) for ${\goth p}'_{n-i}$. When passing to 
the expression in 
the second and third lines
we substitute iterative formulas (\ref{Newton-a}) for the sums with ${\goth p}_{n-i}$ and 
${\goth p}_{n-i-2}$, and 
(\ref{mod-N}) for the sum with ${\goth p}'_{n-i-2}$ (thus using the induction assumption), and
 eq.(\ref{popo}) for 
$(q^2 {\goth p}'_0 - {\goth p}_0)$. We collect 
terms coming from the sums with ${\goth p}_{n-i}$ and ${\goth p}_{n-i-2}$ and containing 
${\goth e}_{n-2i}$, $i\geq 1$, in 
the second line,
and all other terms in the third line. Most of terms in two sums in the second line cancel each other, 
except the term containing ${\goth e}_{n-2}$. In line four we collect all the terms from the second and third 
lines containing ${\goth e}_{n-2}$
and finally observe that they are also cancelling.

\medskip
The right equality in (\ref{mod-N}) is proved analogously.
\hfill$\blacktriangleright$

\medskip
We now notice that in all three cases, considered in theorem \ref{corollary6.4}, the images of the 
elements
${\goth e}_i$ are given by the elementary symmetric functions (see eqs.(\ref{rep-charO+})-(\ref{rep-charOodd})). Hence, by the Wronski relations
(\ref{mod-W}), the images of the modified elements ${\goth h}'_n$ are the complete symmetric
functions in the same arguments. Using then  eq.(\ref{mod-as}) and taking into account  relations
$$
h_i( \nu_0,\{\nu\}_k )= h_{i}(\{\nu\}_k)\,+\,\nu_0\,h_{i-1}(\nu_0,\{\nu\}_k)\ ,
$$
it is easy to check the formulas 
for the images of the elements ${\goth h}_i$ given in eqs.(\ref{para-s-O+}), (\ref{para-s-O-}), (\ref{para-s-O-odd}).
\medskip

To verify formulas for the power sums, we need one more result. 
\begin{lem}  
\lb{lemma6.6} Let ${\goth e}_i$, $i=0,1,2\dots ,$ be a set of commuting independent generators of 
a free commutative algebra.  Assume that another set of the elements ${\goth p}_i$, $i=1,2\dots ,$  
is expressed in terms of the set of the elements ${\goth e}_i$, $i=0,1,2\dots ,$ through the recursion relations 
\be\lb{GS-Newton}
\sum_{i=0}^{n-1}(-q)^i\, {\goth e}_i\, {\goth p}'_{n-i}\, =\, (-1)^{n-1}\, n_q\, {\goth e}_n, \qquad n=1,2,\dots  . 
\ee
Given a natural number $k$, consider the map 
\[ {\goth e}_i \mapsto e_i(\nu_1,\dots ,\nu_k)\ ,\ 0\leq i\leq k\ \ \mbox{and}\ \ {\goth e}_i \mapsto 0\ , 
\ i>k\ ,\] 
to the set of elementary symmetric polynomials in commuting variables $\nu_j$, $j=1,\dots ,k$. 

\medskip
Extending this map by the homomorphicity and
assuming additionally that all differences $(\nu_{r_1}-\nu_{r_2})$, $1\leq r_1<r_2\leq k$, 
are invertible,
we obtain the following expressions for the images of the elements ${\goth p}'_i$
 \be
\lb{GS-expr}
{\goth  p}'_i\, \mapsto\, q^{i-1}\sum_{j=1}^k \delta_j (\nu_j)^i\ ,\ i=1,2,\dots\ ,\quad \mbox{where}\quad \delta_j :=\prod_{r=1\atop r\neq j}^k
{(\nu_j-q^{-2}\nu_r)}/{(\nu_j-\nu_r)}\ .
\ee
Besides, the quantities $\delta_j$ satisfy the relation
\be
\lb{GS-expr2}
q^{-1}\sum_{j=1}^k \delta_j\, =\, q^{-k} k_q.
\ee
\end{lem}
\noindent {\bf Proof.}~
The formula (\ref{GS-expr}) is proved in \cite{GS}, proposition 17. To prove eq.(\ref{GS-expr2}), consider the rational function
in the variable $z$:
\be
\lb{wz1}
w(z):=\prod_{j=1}^k \frac{z-q^{-2}\nu_j}{z-\nu_j}.
\ee
Its expansion in simple ratios is
\be
\lb{wz2}
w(z) = 1\, +\, \sum_{j=1}^k \frac{\alpha_j}{z-\nu_j}, \quad \mbox{where}\quad
\alpha_j=(z-\nu_j)w(z)|_{z=\nu_j}= (1-q^{-2})\, \delta_j \nu_j. 
\ee
Comparing expressions for $w(0)$ following from (\ref{wz1}): $w(0)=\prod_{j=1}^k q^{-2}=q^{-2k}$, and from (\ref{wz2}):
$w(0)=1-(1-q^{-2})\sum_{j=1}^k \delta_j$, one verifies the formula (\ref{GS-expr2}).
\hfill$\blacktriangleright$

\medskip
Now we are ready to derive the parameterization formulas for the power sums. 

\medskip
Let's start with the case $O^+_{2\ell}$.
Lemmas \ref{lemma6.5} and \ref{lemma6.6} fix in this case the spectral parameterization for the modified power sums ${\goth p}'_i$ in the form
\be\lb{deltajo+}
\pi^+_{2\ell}:\;\; {\goth p}'_i\, \mapsto\, q^{i-1}\sum_{j=1}^{2\ell} \delta_j (\nu_j)^i\ , 
\ i\geq 0\ ,\ 
\quad \mbox{where}\quad \delta_j :=\prod_{r=1,\, r\neq j}^{2\ell}
{(\nu_j-q^{-2}\nu_r)}/{(\nu_j-\nu_r)}.
\ee
Here, the formula (\ref{GS-expr2}) confirms the correctness of our choice of the element 
${\goth p}'_0$
 in eq.(\ref{mod-p1}).

\vskip .4cm
Taking the Ansatz
$$
\pi^+_{2\ell}:\;\; {\goth p}_i\, \, \mapsto\, q^{i-1}\sum_{j=1}^{2\ell} d_j (\nu_j)^i\ , \ i\geq 0\ ,\ 
$$
for the images of the power sums ${\goth p}_i$ in ${\cal E}_{ 2\ell}$ we calculate expressions for 
$2\ell$ parameters $d_i$ using eqs.(\ref{mod-p1}). First, rewriting these recursion relations in the 
form~ 
$
{\goth p}_i- g {\goth p}_{i-2} = {\goth p}'_i - q^2 g {\goth p}'_{i-2}
$~
we observe that their images in ${\cal E}_{ 2\ell}$ are fulfilled for any value of index $i$ if we choose
\be
\lb{djo+}
d_j \, =\, {\nu_j^2- \nu_0^2\over \nu_j^2-q^{-2}\nu_0^2}\,\delta_j\, =\, {\nu_j-\nu_{\overline j}\over \nu_j - q^{-2}\nu_{\overline j}}\, \delta_j \, =\, \prod_{r=1\atop r\neq j,{\overline j}}^{2\ell}{\nu_j -q^{-2}\nu_r\over \nu_j-\nu_r}.
\ee
Here we introduced the compact notation $\nu_{\overline j}:=\nu_{2\ell+1-j}$ and used the equality 
$\nu_0^2=\nu_j \nu_{\overline j}$~ and the expression for $\delta_j$ for the transformations. 
Expression for $d_j$ coincides with formula (\ref{para-p-O+}) given in the theorem. To complete a 
verification of the spectral parameterization for the power sums it remains to check the initial 
conditions of the recursion (\ref{mod-p1}). Namely, we have to prove
\be
\lb{spec01-o+}
{\goth p}_0\, \mapsto\, q^{-1}\sum_{j=1}^{2\ell} d_j\, =\, q^{1-2\ell}\Bigl(1+ (2\ell-1)_q\Bigr), \qquad
{\goth p}_1={\goth p}'_1={\goth e}_1\, \mapsto \, \sum_{j=1}^{2\ell} d_j\nu_j\, =\, \sum_{j=1}^{2\ell} \nu_j .
\ee
To this end we consider expansions in the simple ratios of the rational functions $v(z)$ and $v(z)z$: 
$$
v(z):= {\nu_0^2\over z^2 - q^{-2}\nu_0^2} \,u(z), \quad \mbox{where}\quad
u(z) := \prod_{r=1}^{2\ell}\,{z-q^{-2}\nu_r\over z-\nu_r}.
$$

For $v(z)$ one has:~
\ba
\nn
v(z) & =& \displaystyle{\sum_{j=1}^{2\ell} {\alpha_j\over z-\nu_j}\, +\, \beta \Bigl( {1\over z-q^{-1}\nu_0}-{1\over z+q^{-1}\nu_0}\Bigr) },\quad \mbox{where}
\\[2pt]
\nn
\alpha_j & =&  (z-v_j) v(z)|_{z=v_j}=(1-q^{-2})\, \nu_j \delta_j\, {\nu_0^2\over \nu_j^2-q^{-2}\nu_o^2}\, =\,
\nu_j (\delta_j - d_j),
\\[6pt]
\nn
\beta & =&\textstyle\pm (z\mp q^{-1}\nu_0) v(z)|_{z=\pm q^{-1}\nu_0}\, =\,  u(\pm q^{-1}\nu_0)\,\nu_0\, q/2\, =\, \nu_0\,q^{1-2\ell}/ 2.
\ea
Here in the first line notations $\delta_j$ and $d_j$ are defined in (\ref{deltajo+}) and (\ref{djo+}). In 
the second line expressions $u(\pm q^{-1}\nu_0)=q^{-2\ell}$ are derived  using the  conditions 
$\nu_r\nu_{\overline r} = \nu_0^2\;\,\forall\, r=1,\dots, \ell.$ Now evaluating $v(z)$ at $z=0$ in two 
ways (using its multiplicative and additive expressions) we obtain
$
\textstyle
\,-q^{2-4\ell}= v(0) =\sum_{j=1}^{2\ell} (d_j-\delta_j)-q^{2-2\ell}\,
$
which in combination with eq.(\ref{GS-expr2}) proves spectral parameterization formula (\ref{spec01-o+}) for ${\goth p}_0$.

\medskip
In a similar way for $v(z) z$ one calculates: 
$$
v(z) z = \sum_{j=1}^{2\ell}\ {\nu_j^2(\delta_j-d_j)\over z-\nu_j}\, +\, {\nu_0^2 q^{-2\ell}\over 2} 
\Bigl( {1\over z-q^{-1}\nu_0}+{1\over z+q^{-1}\nu_0}\Bigr).
$$
Evaluating this expression at $z=0$ one finds $0 = v(z) z|_{z=0} = \sum_{j=1}^{2\ell} \nu_j (d_j-\delta_j)$ which proves the spectral parameterization formula (\ref{spec01-o+}) for ${\goth p}_1$.
\medskip

Proceeding to the case $O_{2\ell -1}$ we notice that in this case Lemmas \ref{lemma6.5} and 
\ref{lemma6.6} prescribe the following parameterization for the modified power sums:
\be
\lb{deltajo-odd}
\pi_{2\ell-1}:\;\; {\goth p}'_i\, \mapsto\, q^{i-1}\sum_{j=0}^{2\ell-2} \delta_j (\nu_j)^i 
\ , \ i\geq 0\ ,\ 
\quad \mbox{where}\quad \delta_j :=\prod_{r=0,\, r\neq j}^{2\ell-2}
{(\nu_j-q^{-2}\nu_r)}/{(\nu_j-\nu_r)}.
\ee
Taking the ansatz
$$
\pi_{2\ell-1}:\;\; {\goth p}_i\, \, \mapsto\, q^{i-1}\sum_{j=1}^{2\ell-2} d_j (\nu_j)^i + d_0 (\nu_0)^i
\ , \ i\geq 0\ ,\ 
$$
we derive expressions for the parameters~ $d_j, j=1,\dots , 2\ell -2$~ from recursion relations 
(\ref{mod-p1}):
\be
\lb{djo-odd}
d_j \, =\, {\nu_j^2- \nu_0^2\over \nu_j^2-q^{-2}\nu_0^2}\,\delta_j\, =\, {\nu_j-\nu_{\overline j}\over \nu_j - q^{-2}\nu_{\overline j}}\, \delta_j \, =\, {\nu_j-q^{-2}\nu_0\over \nu_j-\nu_0}\prod_{r=1\atop r\neq j,{\overline j}}^{2\ell-2}{\nu_j -q^{-2}\nu_r\over \nu_j-\nu_r}.
\ee
Here $\nu_{\overline j}:=\nu_{2\ell-1-j}$ and we used the relation 
$\nu_0^2=\nu_j \nu_{\overline j}$ to simplify expression for $d_j$. Notice that the parameter $d_0$ is not 
fixed from the recursion relation. It is defined by the initial conditions, which in this case look like
\be
\lb{spec01-odd}
{\goth p}_0\, \mapsto\, q^{-1}\sum_{j=1}^{2\ell-2} d_j + d_0\, =\, q^{2-2\ell}\Bigl(1+ (2\ell-2)_q\Bigr), \quad
\, {\goth p}_1\mapsto \, \sum_{j=1}^{2\ell-2} d_j\nu_j+ d_0\nu_0\, =\, \sum_{j=0}^{2\ell-2} \nu_j .
\ee
To verify these conditions we  expand in the simple ratios functions $v(z)$ and $v(z)z$, where 
$$
v(z):= {z\, \nu_0\over z^2 - q^{-2}\nu_0^2} \, \prod_{r=1}^{2\ell-2}\,{z-q^{-2}\nu_r\over z-\nu_r}.
$$
One has
\ba
\nn
v(z) &=&
\sum_{j=1}^{2\ell-2} {\nu_j(d_j-\delta_j)\over z-\nu_j}\, +\, {\nu_0 q^{2-2\ell}\over 2} \Bigl( {1\over z-q^{-1}\nu_0}+{1\over z+q^{-1}\nu_0}\Bigr),
\\[10pt]
\nn
v(z) z &=& 
\sum_{j=1}^{2\ell-2} {\nu_j^2(d_j-\delta_j)\over z-\nu_j}\, +\, {\nu_0^2 q^{1-2\ell}\over 2} \Bigl( {1\over z-q^{-1}\nu_0}-{1\over z+q^{-1}\nu_0}\Bigr),
\ea
where the expressions $\delta_j$ and $d_j$ are defined in (\ref{deltajo-odd}) and (\ref{djo-odd}). Evaluating 
these functions at $z=0$ and taking into account that, by eq.(\ref{GS-expr2}), 
$\displaystyle{\sum_{j=1}^{2\ell-2}\delta_j = q^{3-2\ell}(2\ell-2)_q}$
we conclude that the equalities in (\ref{spec01-odd}) become true with the choice $d_0= q^{2-2\ell}$. 
Thus we obtained formulas  (\ref{para-s-O-odd}), (\ref{para-p-O-odd}) for $d_j, j=0,\dots, 2\ell-1,$ 
given in the theorem for the case $O(2\ell -1)$.
\medskip

 In the case  $O^-_{2\ell}$ the spectral parameterization of the power sums is suitably derived with the 
 use of the other set $\{{\goth p}''_i\}_{i>0}$ of the 
 modified power sums. First, we notice that for the pair 
 of sets $\{{\goth h}_i\}_{i\geq 0}$ and $\{{\goth p}''_i\}_{i>0}$ satisfying  
the relations (\ref{mod-N}), an analogue of lemma \ref{lemma6.6} is valid. Namely, if one 
 considers the homomorphic map which sends  
the set of the complete sums to the complete symmetric polynomials: 
\[ {\goth h}_i \mapsto h_i(\nu_1, \dots ,\nu_k)\ ,\] 
then the images of the elements ${\goth p}''_i$ are given 
 by the formulas (\ref{GS-expr}), (\ref{GS-expr2}) with 
the substitutuion ${\goth p}'_i\mapsto {\goth p}''_i$. Hence, in the case  $O^-_{2\ell}$ the spectral 
 parameterization for 
 modified power sums ${\goth p}''_i$ looks like 
 \be\lb{deltajo-o-}
 \pi^-_{2\ell}:\;\; {\goth p}''_i\, \mapsto\, q^{i-1}\sum_{j=1}^{2\ell-2} \delta_j (\nu_j)^i\ ,\ i\geq 0\ ,\ 
  \quad \mbox{where}\quad \delta_j :=\prod_{r=1,\, r\neq j}^{2\ell-2}
 {(\nu_j-q^{-2}\nu_r)}/{(\nu_j-\nu_r)}.
 \ee
 For the power sums we take following ansatz
 $$
 \pi^-_{2\ell}:\;\; {\goth p}_i\, \, \mapsto\, q^{i-1}\sum_{j=1}^{2\ell-2} d_j (\nu_j)^i+d_{0,i} (\nu_0)^i 
 \ ,\ i\geq 0\ .
 $$
Here we admit a dependence of the parameter $d_{0,i}$ on a value of index  $i$ of the power sum. 
Recursion relations (\ref{mod-p1}) prescribe following equalities between parameters $d_j$, $d_{0,i}$  and $\delta_j$ in (\ref{deltajo-o-}):
\be
\lb{djo-o-}
d_j \, =\, {\nu_j^2- q^{-4}\nu_0^2\over \nu_j^2-q^{-2}\nu_0^2}\,\delta_j\, =\, {\nu_j-q^{-4}\nu_{\overline j}\over \nu_j - q^{-2}\nu_{\overline j}}\, \delta_j \, =\,  {\nu_j-q^{-4}\nu_{\overline j}\over \nu_j - \nu_{\overline j}}\,\prod_{r=1\atop r\neq j,{\overline j}}^{2\ell-2}{\nu_j -q^{-2}\nu_r\over \nu_j-\nu_r}, \qquad d_{0,i}=d_{0,i-2}.
\ee
Here we denote $\nu_{\overline j}:=\nu_{2\ell-1-j}$. Notice that the conditions above fix values of all 
the parameters $d_{0,i}$ in terms of $d_{0,0}$ and $d_{0,1}$. The latter pair is  fixed by the initial 
conditions: 
\be
\lb{spec01-o-}
{\goth p}_0\, \mapsto\, q^{-1}\sum_{j=1}^{2\ell-2} d_j + d_{0,0}\, =\, q^{1-2\ell}\Bigl(1+ (2\ell-1)_q\Bigr), \quad
{\goth p}_1
\, \mapsto \, \sum_{j=1}^{2\ell-2} d_j\nu_j+ d_{0,1}\nu_0\, =\, \sum_{j=1}^{2\ell-2} \nu_j .
\ee
To verify these conditions we  expand in the simple ratios functions $v(z)$ and $v(z)z$, where now
$$
v(z):= { \nu_0^2\over z^2 - q^{-2}\nu_0^2} \, \prod_{r=1}^{2\ell-2}\,{z-q^{-2}\nu_r\over z-\nu_r}.
$$
One has
\ba
\nn
v(z) &=&
\sum_{j=1}^{2\ell-2} {q^2\nu_j(d_j-\delta_j)\over z-\nu_j}\, +\, {\nu_0 q^{3-2\ell}\over 2} \Biggl( {1\over z-q^{-1}\nu_0}-{1\over z+q^{-1}\nu_0}\Biggr),
\\[10pt]
\nn
v(z) z &=& 
\sum_{j=1}^{2\ell-2} {q^2\nu_j^2(d_j-\delta_j)\over z-\nu_j}\, +\, {\nu_0^2 q^{2-2\ell}\over 2} \Biggl( {1\over z-q^{-1}\nu_0}+{1\over z+q^{-1}\nu_0}\Biggr),
\ea
where the expressions $\delta_j$ and $d_j$ are defined 
in (\ref{deltajo-o-}) and (\ref{djo-o-}). Evaluating 
these functions at $z=0$ and taking into account that, by 
eq.(\ref{GS-expr2}), $\displaystyle{\sum_{j=1}^{2\ell-2}\delta_j = q^{3-2\ell}(2\ell-2)_q}$
we find $d_{0,0}=2 q^{1-2\ell}$, $d_{0,1}=0$ and hence,
$$
d_{0,i} = (1+(-1)^i)q^{1-2\ell}.
$$
Thus we confirmed the spectral parameterization 
formulas  (\ref{para-s-O-}), (\ref{para-p-O-}) given in 
the theorem for the case $O^-_{2\ell}$.
\hfill$\blacksquare$
 
\section*{Acknowledgments}

The authors thank Dimitry Gurevich and Pavel Saponov for their time, fruitful 
discussions and valuable remarks.
The work of the second author (P. P.) was  supported
by the Basic Research Program of HSE University.

\end{document}